\documentclass{gtart_h}

\def\ifplaintex{\expandafter\ifx\csname documentclass\endcsname\relax}

\def\gtp{{\mathsurround=0pt\it $\cal G\mskip-2mu$eometry \&\ 
$\cal T\!\!$opology $\cal P\!$ublications}}  

\def\recd{{\small Received:\qua\receiveddate\ifx\reviseddate\relax
\else\qquad Revised:\qua\reviseddate\fi\par}} 

\def\lognumber#1{\def\thelognumber{#1}}
\def\volumenumber#1{\def\thevolumenumber{#1}}
\def\volumeyear#1{\def\thevolumeyear{#1}}
\def\papernumber#1{\def\thepapernumber{#1}}
\def\pagenumbers#1#2{\def\startpage{#1}\def\finishpage{#2}}
\def\published#1{\def\publishdate{#1}}

\def\received#1{\def\receiveddate{#1}}

\def\accepted#1{\def\accepteddate{#1}}

\def\asciiaddress#1{\def\theasciiaddress{#1}}
\def\asciiemail#1{\def\theasciiemail{#1}}

\long\def\asciiabstract#1{\long\def\theasciiabstract{#1}}

\let\\\par\let\thelognumber\relax\let\thevolumenumber\relax
\let\thepapernumber\relax\let\thevolumeyear\relax\let\startpage\relax
\let\finishpage\relax\let\publishdate\relax\let\receiveddate\relax
\let\reviseddate\relax\let\accepteddate\relax\let\theasciititle\relax
\let\theasciiauthors\relax\let\theasciiaddress\relax
\let\theasciiabstract\relax

\let\theasciiemail\relax

\ifplaintex
\font\logobig=cmssbx10 scaled 3836
\font\logomed=cmssbx10 scaled 2557
\else
\font\logobig=cmssbx10 scaled 4200
\font\logomed=cmssbx10 scaled 2800
\fi

\long\def\makeagttitle{   
\count0=\startpage
\agt\hfill      
\hbox to 45truept{\vbox to 0pt{\vglue -13truept{\logomed A\kern -.37em{\logobig 
T}\kern -.38em G}\vss}\hss}
\break
{\small Volume \thevolumenumber\ (\thevolumeyear)
\startpage--\finishpage\nl
Published: \publishdate}

\vglue .25truein

{\parskip=0pt\leftskip 0pt plus
1fil\def\\{\par\smallskip}{\Large\bf\thetitle}\par\medskip} \vglue
0.05truein

{\parskip=0pt\leftskip 0pt plus 1fil\def\\{\par}{\sc\theauthors}
\par\medskip}%
 
\vglue 0.03truein 

{\small\leftskip 25truept\rightskip 25truept{\bf Abstract}\stdspace\theabstract

{\bf AMS Classification}\stdspace\theprimaryclass
\ifx\thesecondaryclass\relax\else; \thesecondaryclass\fi\par
{\bf Keywords}\stdspace \thekeywords\par}\vglue 7truept

}   

\ifplaintex
\font\phead=cmsl9 scaled 950
\font\pnum=cmbx10 scaled 913
\font\pfoot=cmsl9 scaled 950
\headline{\vbox to 0pt{\vskip -4.5mm\line{\small\phead\ifnum
\count0=\startpage ISSN 1472-2739 (on-line) 1472-2747 (printed)
\hfill {\pnum\folio}\else\ifodd\count0\def\\{ }%
\ifx\theshorttitle\relax\thetitle\else\theshorttitle\fi\hfill{\pnum\folio}
\else\def\\{ and }{\pnum\folio}\hfill\ifx\theshortauthors\relax\theauthors
\else\theshortauthors\fi\fi\fi}\vss}}
\footline{\vbox to 0pt{\vglue 0mm\line{\small\pfoot\ifnum\count0=\startpage
\copyright\ \gtp\hfill\else
\agt, Volume \thevolumenumber\ (\thevolumeyear)\hfill\fi}\vss}}
\else
\font=cmsl9 scaled 1050
\font\lnum=cmbx10 
\font=cmsl9 scaled 1050
\makeatletter
\def\@oddhead{{\small\lhead\ifnum\count0=\startpage ISSN 1472-2739 
(on-line) 1472-2747 (printed)\hfill {\lnum\number\count0}\else\ifodd\count0
\def\\{ }\ifx\theshorttitle\relax \thetitle \else\theshorttitle\fi\hfill
{\lnum\number\count0}\else\def\\{ and }{\lnum\number\count0}
\hfill\ifx\theshortauthors\relax 
\theauthors\else\theshortauthors\fi\fi\fi}}\def\@evenhead{\@oddhead}
\def\@oddfoot{\small\lfoot\ifnum\count0=\startpage\copyright\ \gtp\hfill\else
\agt, Volume \thevolumenumber\ (\thevolumeyear)\hfill\fi}
\def\@evenfoot{\@oddfoot}
\makeatother
\fi
\let\maketitlepage\makeagttitle
\let\makeshorttitle\maketitlepage
\let\maketitle\maketitlepage

\newwrite\gtoutfile
\long\gdef\makeheadfile{  
{\def\\{, }\def\s{ }
\immediate\openout\gtoutfile head.xxx
\immediate\write\gtoutfile{Proxy-for: \ifx\theasciiauthors\relax
\theauthors\else\theasciiauthors\fi\s<\ifx\theasciiemail\relax\theemail\else\theasciiemail\fi>}
\immediate\write\gtoutfile{\noexpand\\}
\immediate\write\gtoutfile{Authors: \ifx\theasciiauthors\relax
\theauthors\else\theasciiauthors\fi}
{\def\\{ }\immediate\write\gtoutfile{Title: \ifx\theasciititle\relax
\thetitle\else\theasciititle\fi}}
\immediate\write\gtoutfile{Subj-class: GT or SG, GR etc}
\immediate\write\gtoutfile{MSC-class: \theprimaryclass\ifx\thesecondaryclass\relax\else, \thesecondaryclass\fi}
\immediate\write\gtoutfile{Journal-ref: Algebr. Geom. Topol. \thevolumenumber\s
(\thevolumeyear) \startpage-\finishpage}
\immediate\write\gtoutfile{Comments: Published by Algebraic and
Geometric Topology at}
\immediate\write\gtoutfile{\s\s\s  http://www.maths.warwick.ac.uk/agt/AGTVol\thevolumenumber/agt-\thevolumenumber-\thepapernumber.abs.html}
\immediate\write\gtoutfile{\noexpand\\}
\immediate\write\gtoutfile{}
\ifx\theasciiabstract\relax
\immediate\write\gtoutfile{\theabstract}\else
\immediate\write\gtoutfile{\theasciiabstract}\fi
\immediate\write\gtoutfile{}
\immediate\write\gtoutfile{\noexpand\\}
\immediate\write\gtoutfile{}
\immediate\closeout\gtoutfile}}  

\def\maketitlepage{\makeagttitle\makeheadfile}
\let\makeshorttitle\maketitlepage
\let\maketitle\maketitlepage

\lognumber{31}
\volumenumber{4}
\volumeyear{2004}
\papernumber{31}
\pagenumbers{685}{719}
\received{6 July 2004} 
\accepted{16 August 2004}
\published{9 September 2004}

\usepackage{graphicx,amssymb,amsmath,amscd}

\newtheorem{theorem}{Theorem}[section]
\newtheorem{proposition}[theorem]{Proposition}
\newtheorem{lemma}[theorem]{Lemma}
\newtheorem{corollary}[theorem]{Corollary}
\theoremstyle{definition}

\newtheorem{definition}[theorem]{Definition}
 
\newtheorem{remark}[theorem]{Remark} 
\newcommand{\zee}{\mbox{$\mathbb Z$}}
\newcommand{\spinc}{\mbox{$Spin^c$ }} 
\newcommand{\s}{\mathfrak s}
\newcommand{\T}{\mathfrak t}
\newcommand{\R}{\mathfrak r}
\newcommand{\sym}{{\mbox{$Sym$}}}
\newcommand{\x}{{\mbox{\bf x}}}

\newcommand{\spin}{$spin^c$-structure }
\newcommand{\spins}{$spin^c$-structures }
\newcommand{\gr}{\mbox{gr} }
\newcommand{\OnS}{{Ozsv\'ath and Szab\'o }}
\newcommand{\mb}{\boldmath}
\newcommand{\cue}{\mbox{$\mathbb Q$}}

\newcommand{\balpha}{\mbox{\mb $\alpha$}}
\newcommand{\bbeta}{\mbox{\mb $\beta$}}
\newcommand{\sbar}{\underline{\s}}
\begin{document}
\title{Heegaard Floer homology of certain mapping tori}
\authors{Stanislav Jabuka\\Thomas Mark}
\addresses{Department of Mathematics, Columbia University\\2990 Broadway, New York, NY 10027, USA\\{\rm and}\\
Department of Mathematics, Southeastern Louisiana University\\1205 North Oak Street, Hammond, LA 70402, USA}
\asciiaddress{Department of Mathematics, Columbia University\\2990 Broadway, New York, NY 10027, USA\\and\\
Department of Mathematics, Southeastern Louisiana University\\1205 North Oak Street, Hammond, LA 70402, USA}
\gtemail{\mailto{jabuka@math.columbia.edu}\qua {\rm and}\qua 
\mailto{Thomas.Mark@selu.edu}}
\asciiemail{jabuka@math.columbia.edu, Thomas.Mark@selu.edu}
\begin{abstract}
We calculate the Heegaard Floer homologies $HF^+(M,\s)$ for mapping tori $M$ 
associated to certain surface diffeomorphisms, where $\s$ is any \spinc 
structure on $M$ whose first Chern class is non-torsion. Let $\gamma$ and 
$\delta$ be a pair of geometrically dual nonseparating curves on a 
genus $g$ Riemann surface $\Sigma_g$, and let $\sigma$ be a curve separating 
$\Sigma_g$ into components of genus $1$ and $g-1$. Write $t_\gamma$, 
$t_\delta$, and $t_\sigma$ for the right-handed Dehn twists 
about each of these curves. The examples we consider are 
the mapping tori of the diffeomorphisms $t_\gamma^m\circ t_\delta^n$ 
for $m,n\in \zee$ and that of $t_\sigma^{\pm 1}$.
\end{abstract}
\asciiabstract{%
We calculate the Heegaard Floer homologies$HF^+(M,s) for mapping tori
M associated to certain surface diffeomorphisms, where s is any Spin^c
structure on M whose first Chern class is non-torsion. Let gamma and
delta be a pair of geometrically dual nonseparating curves on a genus
g Riemann surface Sigma_g, and let sigma be a curve separating Sigma_g
into components of genus 1 and g-1.  Write t-gamma, t_delta, and
t_sigma for the right-handed Dehn twists about each of these curves.
The examples we consider are the mapping tori of the diffeomorphisms
t_gamma^m circ t_delta^n for m,n in Z and that of t_sigma^{+-1}.}
\primaryclass{57R58}
\secondaryclass{53D40}
\keywords{Heegaard Floer homology, mapping tori}
\maketitlepage

\section{Introduction} 

In \cite{peter4}, Peter Ozsv\'ath and Zolt\'an Szab\'o introduced a
set of new invariants of 3-manifolds, the Heegaard Floer homology
groups. There are several variations available in the construction,
which give rise to several related invariants $HF^+(Y,\s)$,
$HF^-(Y,\s)$, $HF^\infty(Y,\s)$, and $\widehat{HF}(Y,\s)$. Each of
these is a relatively graded group associated to a closed oriented
3-manifold $Y$ equipped with a \spinc structure $\s$. This paper is
concerned with the calculation of the group $HF^+$ in case $Y$ is a
fibered 3-manifold whose monodromy is of a particular form. 

If $\Sigma_g$ is a closed oriented surface of genus $g>1$ and
$\phi:\Sigma_g\to \Sigma_g$ is an orientation-preserving
diffeomorphism, we can form the mapping torus $M(\phi)$ as a quotient
of $\Sigma_g\times[0,1]$. For a simple closed curve $c$ on $\Sigma_g$,
let $t_c:\Sigma_g\to \Sigma_g$ denote the right-handed Dehn twist
about $c$. Choose a pair of geometrically dual nonseparating circles
$\gamma$ and $\delta$ on $\Sigma_g$, and also let $\sigma$ be a circle
that separates $\Sigma_g$ into components of genus $1$ and $g-1$. Our
main result is the calculation of $HF^+(M, \s)$, where $M$ is any of
the mapping tori $M(t_\gamma^m t_\delta^n)$  for
$m,n\in\zee$ or $M(t_\sigma^{\pm 1} )$  and $\s$ is any non-torsion \spinc structure on $M$.

To state the results, recall that the homology group
$H_2(M(\phi);\zee)$ of the mapping torus $M(\phi)$ can be identified
with $\zee\oplus \ker(1-\phi_*)$ where $\phi_*$ denotes the action of
$\phi$ on $H_1(\Sigma;\zee)$. For a fixed integer $k$, the two
requirements
\begin{enumerate}
\item $\langle c_1(\mathfrak{s}) , [\Sigma _g]\rangle = 2k$ and
\item $\langle c_1(\mathfrak{s}) , [T]\rangle = 0$ for all classes
$[T]$ coming from $H_1(\Sigma_g)$
\end{enumerate}
specify the \spinc structure $\s$ uniquely modulo torsion. Let ${\mathcal S}_k\subset 
Spin^c(M(\phi))$ denote the collection of \spinc structures 
satisfying these two conditions: then for $\s\in{\mathcal S}_k$ the adjunction
inequality for Heegaard Floer homology shows that $HF^+(M(\phi),\s) =
0$ unless $|k|\leq g-1$.

Now let $X(g,d)$ denote the $\zee$-graded group 
whose value in degree $j$ is the homology group
$H_{g-j}(\sym^d(\Sigma_g);\zee)$ of the $d$-th symmetric power of 
$\Sigma_g$. By
convention, $X(g,d) = 0$ if $d<0$. We denote by $X(g,d)[m]$ the
graded group $X(g,d)$ with the grading shifted up by $m$, and if $G$
is a group we let $G_{(m)}$ denote the graded group $G$ concentrated
in grading $m$.

\begin{theorem}\label{mainmain}
Fix an integer $k$ with $0<|k|\leq g-1$, and let $n$ be a nonzero 
integer. Define $d = g-1-|k|$. For the mapping torus $M(t_\gamma^n)$ of the composition of 
$n$ right-handed Dehn twists about a nonseparating circle 
$\gamma\subset \Sigma_g$, there is a unique \spinc structure 
$\s_k\in{\mathcal S}_k$ for which we have an isomorphism of relatively 
graded groups
\[
HF^+(M(t_\gamma^n),\s_k) \cong (X(g-1,d-1)\otimes 
H^*(S^1;\zee))[\varepsilon(n)] \oplus 
\Lambda^{2g-2-d}H^1(\Sigma_{g-1})_{(g-d)}.
\]
Here $\varepsilon(n) = 0$ if $n>0$ and $\varepsilon(n) = -1$ if $n<0$. 

For any other $\s'\in{\mathcal S}_k$, we have an isomorphism
\[
HF^+(M(t_\gamma^n),\s') \cong X(g-1,d-1)\otimes H^*(S^1;\zee).
\]

\end{theorem}

\begin{theorem}\label{mainmain2}
For $k$ and $d$ as above, let $m,n$ be nonzero integers. Then for the 
mapping torus $M(t_\gamma^mt_\delta^n)$ there is 
a unique \spinc structure $\s_k\in{\mathcal S}_k$ for which there is an 
isomorphism of relatively graded groups
\[
HF^+(M(t_\gamma^mt_\delta^n), \s_k) \cong\left\{
\begin{array}{l}
X(g-1,d-1)\oplus \Lambda^{2g-2-d}H^1(\Sigma_{g-1})_{(g-1-d)} \\
X(g-1,d-1)[-1]\oplus\Lambda^{2g-2-d}H^1(\Sigma_{g-1})_{(g-1-d)}\\
X(g-1,d-1)[-2]\oplus\Lambda^{2g-2-d}H^1(\Sigma_{g-1})_{(g-1-d)} 
\end{array}\right.
\]
The three lines on the right-hand side above correspond to the three cases 
$m,n>0$, $m\cdot n <0$ and $m,n<0$ respectively. 
For all other $\s'\in{\mathcal S}_k$ (and regardless of the signs of $m$ 
and $n$) we have an isomorphism
\[
HF^+(M(t_\gamma^mt_\delta^n),\s') \cong X(g-1,d-1).
\]
\end{theorem}

For the case of the separating curve $\sigma\subset\Sigma_g$, there 
is no torsion in the second cohomology. Thus there is a unique \spinc 
structure $\s_k\in{\mathcal S}_k$, and we have

\begin{theorem}\label{mainmain3}
Let $n$ be either $1$ or $-1$  and let $k$ and $d$ be as above. Then for 
a genus-1 separating curve $\sigma$, we have an isomorphism of 
relatively graded groups
\[
HF^+(M(t_\sigma^n) , \mathfrak{s} _k) \cong (X(g-1,d-1) \otimes
H^*(S^1 \sqcup S^1 )){[\varepsilon(n)]} \oplus
\Lambda^{2g-2-d}H^1(\Sigma_{g-1})_{(g-d)},
\]
where $\varepsilon(n)$ is as in theorem \ref{mainmain}.
\end{theorem}

\begin{remark} (1)\qua When $\s$ is not a torsion \spinc structure, the
Heegaard Floer homology group $HF^+(M,\s)$ carries only a relative
cyclic grading. Thus the above results should be interpreted as
asserting the existence of a $\zee$ grading on the relevant Floer
homology groups such that the stated isomorphisms hold.

(2)\qua If $\s\notin {\mathcal S}_k$ for some $k$, then for $\phi$ one of the
diffeomorphisms considered above we have $HF^+(M(\phi),\s) = 0$ by the
adjunction inequality: the homology of $M(\phi)$ is spanned by
$\Sigma_g$ and classes represented by tori. The above theorems
therefore give a complete calculation of $HF^+(M(\phi))$ in nontorsion
\spinc structures in the cases listed.

(3)\qua The results of theorems \ref{mainmain}--\ref{mainmain3} are sharpened
slightly in theorems \ref{main2}, \ref{main3a}, \ref{main3b}, and
\ref{main1}. In particular, we determine the ``special'' \spinc
structure $\s_k$ precisely, and in some cases we also describe the
action of $H_1(M;\mathbb{Z}) \otimes _\mathbb{Z} \mathbb{Z}[U]$ on
$HF^+(M;\mathfrak{s})$.
\end{remark}

As an application of theorems \ref{mainmain}--\ref{mainmain3}, we
make the following observation on the genera of the possible fibration
structures on $Y$, where $Y$ is one of the mapping tori considered in
the theorems above. We should remark that the following result can
also be deduced from the relation $\|\cdot\|_A \leq \|\cdot\|_T$
between the Alexander and Thurston norms on $H^1(Y;\zee)$ due to
McMullen \cite{mcmullen}, together with the relationships between the
Euler characteristic of $HF^+$, the Seiberg-Witten invariant, and the
Alexander polynomial. We give a proof based exclusively on Heegaard
Floer homology.

\begin{theorem}\label{fibrthm}
Let $Y$ denote the mapping torus of one of the diffeomorphisms
$t_\gamma^mt_\delta^n$ or $t_\sigma^{\pm 1}$ of the surface $\Sigma_g$ of
genus $g\geq 2$ as above (including the identity diffeomorphism), and
suppose $f:Y\to S^1$ is a fibration of $Y$ having as fiber the
connected surface $F$. Then the genus of $F$ is of the form $g+n(g-1)$
for some $n\geq 0$.
\end{theorem}

\begin{proof}
In \cite{peter6}, Ozsv\'ath and Szab\'o show that if $Y$ is a
3-manifold that fibers over the circle with connected fiber $F$ of genus
$h>1$, then there is a unique \spinc structure $\s$ on $Y$ satisfying
the following conditions:
\begin{enumerate}
\item $\langle c_1(\s), [F]\rangle = 2h-2,$ and
\item $HF^+(Y, \s)\neq 0$.
\end{enumerate}
Furthermore, they prove that for this ``canonical'' \spinc structure,
$HF^+(Y,\s)\cong \zee$. Theorem \ref{fibrthm} follows from this
together with the observation that according to theorems
\ref{mainmain}--\ref{mainmain3}, there is only one possible canonical
\spinc structure. To see this, note that if $Y$ is $M(t_\gamma^n)$ 
for some $n\in \zee$ or
$M(t_\sigma^{\pm 1})$  then there is only one \spinc
structure (up to conjugation) having $HF^+$ equal to $\zee$, namely
$\s_{g-1}$. If $Y$ is $M(t_\gamma^mt_\delta^n)$ there may be more such
\spinc structures, but note that since all elements of ${\mathcal S}_k$
differ by torsion classes, all the Chern classes of \spinc structures
in that set have the same pairing with every homology class $[F]\in
H_2(Y;\zee)$. Hence the uniqueness of the \spinc structure satisfying
(1) and (2) above implies that if ${\mathcal S}_k$ contains a canonical
\spinc structure then it contains only that \spinc structure. This
with an examination of theorem \ref{mainmain2} rules out all $k$
except $k = g-1$, so $\s_{g-1}$ is the only possible canonical \spinc
structure.

Now suppose $f:Y\to S^1$ is a fibration with connected fiber $F$,
where $Y = \Sigma_g\times [0,1]/\sim$ is a mapping torus as in the
statement. The canonical \spinc structure $\s_{g-1}$ on $Y$ has
$\langle c_1(\s_{g-1}),[\Sigma_g]\rangle = 2g-2$, and $\langle
c_1(\s_{g-1}), [T]\rangle = 0$ for all classes $[T]\in H_2(Y; \zee)$
coming from $\ker(1-\phi_*)$. It follows that the image of
$c_1(\s_{g-1})$ under the natural homomorphism
\[
\rho: H^2(Y;\zee)\to \mbox{Hom}(H_2(Y;\zee), \zee)
\]
is divisible by $2g-2$. From the remarks above, we know $\s_{g-1}$
is also the canonical \spinc structure for the
fibration $f$: in particular we must have $\langle c_1(s_{g-1}),
[F]\rangle = 2h-2$, where $h$ is the genus of $F$. Hence $2h-2$ is
divisible by $2g-2$, which is equivalent to the statement that $h = g+
n(g-1)$ for some $n\geq 0$.
\end{proof}

\begin{remark} We implicitly assume here that $h\geq 2$. However, if a
3-manifold $Y$ admits a fibration with fiber genus 1---ie, $Y$ is a
torus bundle over $S^1$---then $H_2(Y;\zee)$ is generated by classes
represented by tori. It then follows from the adjunction inequality
that any \spinc structure having nontrivial $HF^+$ must be torsion, so
from theorems \ref{mainmain}--\ref{mainmain3} none of the examples
under consideration can admit such a structure.
\end{remark}

It is not hard to construct fibrations of genus $g + n(g-1)$ on $Y$ as
above for any $n\geq 0$: see for example \cite{tollefson},
\cite{neumann}. In fact, for any class $[F]$ with $\langle
c_1(\s_{g-1}),[F]\rangle \neq 0$ there is a fibration of $Y$ having
fiber with homology class $[F]$ (see \cite{thurston}).

It is interesting to compare the results of theorems
\ref{mainmain}--\ref{mainmain3} to Floer homologies of different flavours
(see for example \cite{eaman}, \cite{hutchings}, \cite{seidel}).
Specifically, Eaman Eftekhary \cite{eaman} and Paul Seidel
\cite{seidel} have calculated the Lagrangian Floer homology $HF(f)$ of
a surface diffeomorphism $f$ obtained by Dehn twists along an
arrangement of curves $C=C_1^+\cup \ldots \cup C_r^+\cup C_1^-\cup \ldots
\cup C_s^- \subset \Sigma _g$ (the superscripts indicate the
handedness of the Dehn twists performed). Their results show that
\[
HF(f) \cong  H^* (\Sigma _g \backslash (\cup _i C_i^-), \cup _j C_j ^+ )  
\]
Our formulas from theorems \ref{mainmain}--\ref{mainmain3} agree with the
above when $k=g-2$.  

\begin{corollary}
With the notation and hypotheses as in theorems
\ref{mainmain}--\ref{mainmain3}, there
are isomorphisms of relatively graded groups
\begin{align} \nonumber
HF^+(M(t_\gamma^n) , \mathfrak{s}_{g-2}) = & \left\{
\begin{array}{ll}
  H^* (\Sigma _g \backslash \gamma) \quad & n=-1 \cr
 H^* (\Sigma _g,  \gamma ) \quad & n=1
\end{array} \right. \cr 
HF^+(M(t_\gamma^m t_\delta^n) , \mathfrak{s}_{g-2}) = & \left\{
\begin{array}{ll}
H^*(\Sigma _g, \gamma \cup \delta ) \quad & m=n =1 \cr
H^*(\Sigma _g \backslash  \gamma , \delta ) \quad & m\cdot n =- 1 \cr
H^*(\Sigma _g \backslash (\gamma \cup \delta)  ) \quad & m=n =-1 \cr
\end{array}
\right.\cr
HF^+(M(t_\sigma^n) , \mathfrak{s}_{g-2}) = & \left\{
\begin{array}{ll}
  H^* (\Sigma _g \backslash \sigma ) \quad & n=-1 \cr
 H^* (\Sigma _g,  \sigma ) \quad & n=1
\end{array} \right. \cr 
\end{align}
\end{corollary}

The corollary follows from theorems \ref{mainmain}--\ref{mainmain3}
by direct calculation. For example, theorem \ref{mainmain} gives 
$HF^+(M(t_\gamma),\s_{g-2}) \cong \mathbb{Z}_{(g-1)}^{2g-1}\oplus \mathbb{Z}_{(g)}$  which is 
up to degree shift isomorphic to $H^* (\Sigma _g,  \gamma )$. 

Likewise, it is not hard to see using theorem \ref{mainmain} 
that $HF^+(M(t_\gamma),\s_k)$ agrees 
with the periodic Floer homology $HP_*$ of Hutchings and Sullivan 
\cite{hutchings}, in cases where the latter has been calculated. In 
particular, compare theorem \ref{mainmain} above with theorem 5.3 of 
\cite{hutchings}.

It is also worth comparing our results to the Seiberg-Witten
invariants of the mapping tori. Recall that the Euler characteristic $\chi
(HF^+(M,\mathfrak{s}))$ calculates the Seiberg-Witten invariant
$SW_M(\mathfrak{s})$ (this follows from results in \cite{peter3}
together with \cite{mengtaubes}). Theorem \ref{mainmain} 
together with the result
$HF^+(M({\mbox{id}}),\mathfrak{s}_k) \cong X(g,d)$ from \cite{peter1}
(where $k\ne 0$ and $d=g-1-|k|$) give the following observation:

\begin{corollary} \label{onepointeight}
The Seiberg-Witten invariants of the mapping tori $M(\mbox{id})$ and 
$M(t_\gamma^n)$ with $0\neq n\in\zee$ in the
spin$^c$-structure $\mathfrak{s}_k$ are all equal to $\pm {2g-2
\choose g-1-|k|}$ while the Heegaard Floer homology groups in these \spinc
structures are mutually nonisomorphic when $|k|\ne g-1$.
\end{corollary}
Note that the Euler characteristics of the groups $HF^+(M(\phi),\s')$,
for $\phi$ as in the corollary and $\s_k\neq \s'\in{\mathcal S}_k$, are
all zero.
\begin{proof} In section \ref{five5} we show that $HF^+(M(t_\gamma^n),\mathfrak{s}_k)$ is calculated 
as the homology of $(X(g,d),D)$ where $D$ is a certain differential on $X(g,d)$ (which is described in 
detail in section \ref{notationshort}). In section \ref{notationshort} we establish an isomorphism of  
$X(g,d)$ with $X_+(g,d)\oplus X(g-1,d-1)\oplus X(g-1,d-1)$ and show that of the last two summands, only 
the first one consists of cycles for $D$. Thus the total rank of $HF^+(M(t_\gamma^n),\mathfrak{s}_k)$
is less than that of $HF^+(M(\mbox{id}),\mathfrak{s}_k) = X(g,d)$ by at least the rank of 
$X(g-1,d-1)$. The latter has positive rank precisely when $d-1 \ge 0$ or equivalently when $|k| \le g-2$.  
\end{proof}

The main tools used to prove theorems \ref{mainmain}--\ref{mainmain3}
are introduced in \cite{peter1}, where \OnS extend their theory to
produce invariants of nullhomologous knots in 3-manifolds. In the same
article they developed a technique (which we describe in some detail
in section \ref{prelim}) that under favorable circumstances allows one
to calculate the Heegaard Floer homology of a manifold obtained as the
zero framed surgery along a knot $K$ in a 3-manifold $Y$. \OnS used
this technique to find the Floer homology groups of the mapping tori
associated to the identity map and also the diffeomorphism induced by
a single Dehn twist along a non-separating curve (though in the latter
case the result is not presented explicitly). The proofs of the
theorems above follow the outline of those calculations.

The remainder of the paper is organized as follows. In section
\ref{prelim} we review the construction from \cite{peter1} which
explains how to calculate $HF^+(Y_0(K))$ from $HF^+(Y)$ and a certain
quotient complex of the ``knot complex'' $CFK^\infty (Y,K)$. In
section \ref{notationshort} we introduce some auxiliary notation that
is used in subsequent sections and also give a model calculation for
$HF^+(M(t_\gamma))$ that sets the agenda for later sections. In
sections \ref{five5}, \ref{six6} and \ref{four4} we calculate the
Heegaard Floer homologies of $M(t_\gamma^n)$, $M(t_\gamma^m
t_\delta^n)$ and $M(t_\sigma^{\pm 1})$ respectively.

\medskip{\bf Acknowledgements}\qua We would like to thank Peter Ozsv\'ath for
useful conversations and encouragement, and Ron Fintushel for
referring us to \cite{tollefson}.

\section{Calculating the Heegaard Floer homology  of the zero surgery} \label{prelim}
In this section we gather some general results that are used in our
calculations in the subsequent sections. The results presented here
are an adaptation of the techniques developed in \cite{peter1},
tailored for the applications we have in mind. The main theorem of
this section (theorem \ref{morerelevantone} below) explains how one
can calculate the Heegaard Floer homology $HF^+$ of the manifold
$Y_0(K)$ obtained from zero surgery along a nullhomologous knot $K$ in
$Y$. Recall that a knot $K\subset Y$ gives rise to a refined version
of Heegaard Floer homology, whose construction we now summarize. The
reader is referred to \cite{peter1} for more details.

Given $K$, let $E$ denote the torus boundary of a regular neighborhood
of $K$. One can then find a Heegaard surface for $Y$ of the form $E\#
\Sigma_{g-1}$, with attaching circles $\balpha = \alpha_1,\ldots,\alpha_g$
and $\bbeta = \beta_1,\ldots,\beta_g$. Fix a meridian for $K$ lying on
$E$, and let $w$ and $z$ denote a pair of basepoints, one on each side
of this meridian. The data $(E\#\Sigma, \balpha,\bbeta, w)$ together
with a choice of \spinc structure $\s$ on $Y$ can be used to define the
Heegaard Floer chain groups $CF^\infty(Y, \s)$. The additional
basepoint $z$, along with a choice of ``relative \spinc structure''
$\underline{\s}\in Spin^c(Y_0(K))$ lifting $\s$ gives rise to a filtration
$\mathcal F$ on $CF^\infty(Y,\s)$. The ``knot chain complex''
$CFK^\infty(Y,K,\sbar)$ is this filtered complex.

More concretely, we assume that $K$ is nullhomologous and fix a
Seifert surface $F$ for $K$. Then $F$ specifies the zero-framing on
$K$, and can be capped off to a closed surface $\hat{F}$ in the
zero-surgery $Y_0(K)$. The generators of $CFK^\infty(Y,K,\sbar)$ are
triples $[\x,i,j]$, where $\x$ denotes an intersection point between
the tori $T_\alpha = \alpha_1\times\cdots\times
\alpha_g$ and $T_\beta = \beta_1\times \cdots\times\beta_g$ in the
symmetric power $\sym^g(E\#\Sigma)$, and $i$ and $j$ are integers.
The point $\x$ along with the basepoint $w$ determine a \spinc
structure $\s_w(\x)$ on $Y$ as well as a relative \spinc structure
$\underline{\s}_w(\x)$; we require that ${\s}_w(\x) =
{\s}$. Furthermore, $i$ and $j$ are required to satisfy the equation
\[
\langle c_1(\sbar_w(\x)), [\hat{F}]\rangle = 2(j-i).
\]
 In this notation, the filtration
$\mathcal F$ is simply ${\mathcal F}([\x,i,j]) = j$; changing the Seifert
surface $F$ shifts $\mathcal F$ by a constant.

The boundary map $\partial^\infty$ in $CFK^\infty$ is
defined by counting holomorphic disks in $\sym^g(E\#\Sigma)$ and can
only decrease the integers $i$ and $j$. Thus, for example, the subgroup 
$C\{i<0\}$ of $CFK^\infty$ generated by those $[\x,i,j]$ having $i<0$
is a subcomplex, and indeed is simply $CF^-(Y,\s)$ with an additional
filtration. The quotient of $CFK^\infty$ by $C\{i<0\}$ is written
$C\{i\geq 0\}$, and is a filtered version of $CF^+(Y,\s)$. We will use
other similar notations to indicate other sub- or quotient complexes
of $CFK^\infty$. In particular, $\widehat{CFK}$ is by definition the
quotient complex $C\{i = 0\}$, and $\widehat{HFK}$ is the homology of
the graded object associated to the filtration ${\mathcal F}$ of
$\widehat{CFK}$. We denote by $\widehat{HFK}(Y,K;j)$ the summand of
this group supported in filtration level $j$ (typically suppressing
the \spinc structure from the notation).

As an additional piece of structure, we have a natural chain endomorphism
$U$ on $CFK^\infty$ given by $U:[\x,i,j]\mapsto [\x,i-1,j-1]$.

One of the main calculational tools available in Floer homology is the
long exact sequence in homology arising from a triple
of manifolds $Y$, $Y_0$, and $Y_n$. Here $Y_0$ and $Y_n$ denote the
results of $0$- or $n$-framed surgery on a nullhomologous knot $K$ in
$Y$. In this situation, \OnS define a surjective map $Q:Spin^c(Y_0)\to
Spin^c(Y_n)$ (see section 9 of \cite{peter3}), and set
\[
HF^+(Y_0, Q^{-1}(\T)) = \bigoplus_{Q(\s) = \T} HF^+(Y_0, \s).
\]
Furthermore, if we let $W_n$ denote the cobordism connecting $Y_n$
with $Y$ comprising a single 2-handle addition, then a Seifert
surface $F$ for $K$ can be completed to a surface $\tilde{F}$ in
$W_n$. Then given a \spinc structure $\s$ on $Y$ and an integer $k$, we can
specify a \spinc structure $\s_k$ on $Y_n$ by requiring that $\s_k$ be
cobordant to $\s$ by a \spinc structure $\R$ on $W_n$ having $\langle
c_1(\R), [\tilde{F}]\rangle = n+2k$.

The theorem we're aiming for in this section is:

\begin{theorem} \label{morerelevantone}
Let $K$ be a nullhomologous knot in $Y$, and assume that for all
sufficiently large positive integers $n$, $HF^+_{red}(Y) =
HF^+_{red}(Y_n) = 0$ . Fix a positive integer $k$. Then there is a
grading on $HF^+(Y_0 , Q^{-1}(\s_k))$ as a $\Lambda ^* H_1(Y_0)\otimes
_\mathbb{Z} \mathbb{Z}[U]$-module for which there is an isomorphism
$$HF^+(Y_0 , Q^{-1}(\s_k)) \cong H_*(C\{ i<0 \mbox{ and } j\ge k \})$$
of graded $\mathbb{Z}$-modules. The case of $k<0$ is handled by conjugation duality of $HF^+$ (see comment below). 
\end{theorem}

Note that if $\mathfrak{s}_k  = (\mathfrak{s},k) \in \mbox{Spin}^c(Y_n)$ is determined by 
$\mathfrak{s}\in \mbox{Spin}^c(Y)$ and $k\in \mathbb{Z}$, then the dual \spin $\overline{\mathfrak{s}_k}$ is 
given by $\overline{\mathfrak{s}_k} = (\bar{\mathfrak{s}}, -k) = \bar{\mathfrak{s}}_{-k}$ . 
On the other hand, 
the map $Q:\mbox{Spin}^c(Y_0) \rightarrow \mbox{Spin}^c(Y_n)$ is 
equivariant with respect to conjugation, specifically 
$$Q^{-1}(\mathfrak{s}_{-k}) = \overline{ Q^{-1} ( \bar{\mathfrak{s}}_k) } $$
The invariance of $HF^+$ under conjugation of \spins together with theorem \ref{morerelevantone}
now gives 
\begin{equation} \label{negativeconversion}
HF^+(Y_0, Q^{-1}(\mathfrak{s}_{-k})) \cong HF^+(Y_0, Q^{-1}(\bar{\mathfrak{s}}_{k}))
\end{equation}
In what follows we will thus focus on the case $k>0$, the case $k<0$ is handled by equation 
\eqref{negativeconversion}. 

Notice furthermore that by choosing $n$ sufficiently large, 
we can assume that $Q^{-1}(\mathfrak{s}_k)$ consists of a single \spin 
with nonvanishing $HF^+$. 
The reason for this is that under the identification $\mbox{Spin}^c(Y_0) 
\cong \mbox{Spin}^c(Y) \oplus \mathbb{Z}$, any two \spins $\mathfrak{t}_1, \mathfrak{t}_2 
\in Q^{-1}(\mathfrak{s}_k)$ differ by $(0,2n)$. Since $Y_0$ has only finitely many \spin 
$\mathfrak{t}$ for which $HF^+(Y_0,\mathfrak{t}) \ne 0$, the assertion follows. 

We also consider conditions under which the isomorphism from theorem \ref{morerelevantone} is an
isomorphism of $\Lambda^*H_1(Y_0)\otimes \zee[U]$ modules: see theorem
\ref{modulethm} below.

We should note that this theorem was proved in \cite{peter1} for the
case relevant to the mapping tori $M(\mbox{id})$ and $M(t_\gamma)$,
and that the proof given here is little more than a summary of that one.

\subsection{$HF^+(Y_0(K))$ as a $\mathbb{Z}$-module} \label{twoone}

We proceed to the proof of theorem \ref{morerelevantone}.

 As above, a choice
of integer $k$ gives rise to a \spinc structure $\s_k$ on $Y_n$. In
what follows, we will write $CF^+(Y_n, [k])$ for $CF^+(Y_n, \s_k)$,
and often suppress the \spinc structure entirely from the notation for the
knot complex. In section 4 of \cite{peter1}, Ozsv\'ath and Szab\'o prove the
following theorem which calculates $HF^+(Y_n)$ in terms of the
homology of a quotient complex of $CFK^\infty (Y,K)$.
\begin{theorem}[Theorem 4.4 of \cite{peter1}] \label{hello} For all $n$
sufficiently large, there exists an isomorphism of chain complexes
$${^b\Psi^+} : CF^+ (Y_n , [k]) \rightarrow C\{ i \ge 0 \mbox{ or } j \ge k \} $$
\end{theorem}

We proceed by considering the surgery long exact sequence for
$Y\rightarrow Y_0 \rightarrow Y_n$ for $n$ chosen sufficiently large
so that theorem \ref{hello} applies. The sequence in question is
\begin{equation} \label{LES} 
\cdots\stackrel{F}{\longrightarrow} HF^+(Y)
\stackrel{G}{\longrightarrow} HF^+(Y_0 , Q^{-1}([k]))
\stackrel{H}{\longrightarrow} HF^+(Y_n , [k])
\stackrel{F}{\longrightarrow} \cdots
\end{equation}

where $Q:Spin^c(Y_0) \rightarrow Spin^c(Y_n)$ is the surjective map
mentioned previously.

\begin{lemma} \label{auxi}
If $HF^+_{red}(Y) = 0$ and $k> 0$, then  the group 
$HF^+(Y_0,Q^{-1}([k]))$ in \eqref{LES}  is isomorphic to $\ker F$.
\end{lemma}

\begin{proof} 
The assumption $k>0$ implies that the map \[HF^+(Y_0,Q^{-1}([k]))
\rightarrow HF^+_{red}(Y_0,Q^{-1}([k]))\] is an isomorphism (see
Corollary 2.4 in \cite{peter6}). This shows that in the commutative
diagram

\begin{equation} \nonumber
\begin{CD} 
HF^\infty(Y) @>\tilde{G} >> HF^\infty (Y_0 , Q^{-1}([k])) \\
@V{\pi _1 }VV                            @VV{\pi _2}V \\
 HF^+(Y)      @>{G }>>          HF^+(Y_0,Q^{-1}([k])) \\
\end{CD}  
\end{equation}
the map $\pi _2$ is the zero map. On the other hand, our assumption
$HF^+_{red} (Y) = 0$ shows that $\pi _1$ is surjective. Combined,
these two observations imply that the map $G: HF^+(Y)
\to HF^+(Y_0,Q^{-1}([k]))$ is the zero map and so the claim of the lemma follows
from the exactness of \eqref{LES}.
\end{proof}

Now, the map 
\[
F= F_{W_n} : HF^+(Y_n,[k]) \rightarrow HF^+(Y)
\]
appearing in the surgery long exact sequence is the sum of maps 
induced by various {\spinc} structures on the cobordism $W_n$. 
Specifically, we can write $F = f_1 + f_2$, where $f_1$ is the
component of $F$ induced by the spin$^c$-structure $\mathfrak{r}$ on
$W_n$ extending $\s$ and having $\langle c_1(\mathfrak{r} ) ,[\tilde{F}]
\rangle = 2k-n$, and $f_2$ is the sum of the maps induced by
spin$^c$-stuctures $\mathfrak{r}$ with $\langle c_1(\mathfrak{r} )
,[\tilde{F}] \rangle = 2k+(2\ell -1) n$ with $\ell \ne 0$. Let us
assume that the \spinc structure on $Y$ (and therefore the relevant
ones on $Y_n$) is torsion, so that there is an absolute grading on
$HF^+$. This will be the case in all situations we'll consider. Then the
degree shift formula (cf. \cite{peter2}) for maps induced by
cobordisms shows that each component of $f_2$ is homogeneous of degree
at least $2k$ smaller than that of $f_1$ (see the proof of theorem 9.1 in 
\cite{peter1}). Using this fact, Ozsv\'ath
and Szab\'o prove (in section 9 of \cite{peter1}) that there is a
$\mathbb{Z}$-module identification $\ker F \cong \ker f_1$. In their
case the manifold $Y$ is $\#^{2g} (S^1\times S^2)$ but the only
property of this manifold required by their proof is the fact
$HF^+_{red}(Y) = 0$, an assumption that we build into our discussion.
On the other hand, $\ker f_1$ can be identified as the homology of a
quotient complex of $CFK^\infty (Y,K)$. Namely, consider the short
exact sequence
\begin{equation} \label{SES}
 0 \rightarrow C\{ i < 0 \mbox{ and } j\ge k \} \stackrel{\sigma}{\longrightarrow}  C\{ i \ge 0 \mbox{ or } j\ge k\} 
\stackrel{\tau}{\longrightarrow} C\{ i \ge 0\} \rightarrow 0 
\end{equation}
Note that $C\{i\geq 0\}$ is simply a filtered version of $CF^+(Y)$.
With this in mind, we have that $\ker f_1$ can be identified (as a
$\mathbb{Z}$-module) with $\ker \tau_*$. This is also proved in
section 9 of \cite{peter1} for the case $Y=\#^{2g} (S^1\times S^2)$
and it translates verbatim to general case provided one assumes
$HF^+_{red}(Y) = 0$. This latter assumption also implies that the
connecting homomorphism $\delta : H_*(C\{ i \ge 0\}) \rightarrow
H_*(C\{ i < 0 \mbox{ and } j\ge k \})$ is trivial and so $\ker \tau _*
\cong \mbox{Im}(\sigma _*) \cong H_*(C\{ i < 0 \mbox{ and } j\ge k
\})$. These remarks prove theorem \ref{morerelevantone}.

\subsection{The module structure on $HF^+(Y_0(K))$} 
Heegaard Floer homology groups come with an additional algebraic
structure, namely an action of $H_1(M;\mathbb{Z} )/Tors \otimes
_\mathbb{Z} \mathbb{Z} [U]$. In this section we describe this action
in certain cases for $M= Y_0(K)$ (see theorem \ref{modulethm} below).

Observe that there are isomorphisms 
$$\begin{array}{rl}
H_1(Y_0(K)) /Tors & \cong ( H_1(Y)/Tors) \oplus \mathbb{Z} \cr 
H_1(Y_n(K))/ Tors & \cong H_1(Y)/Tors
\end{array}$$
In particular, $H_1(Y)/Tors$ acts on all three terms appearing in
the long exact sequence \eqref{LES}. It follows from the fact that
$F$, $G$, and $H$ in that sequence are induced by cobordisms that they
are all equivariant with respect to the $H_1(Y)/Tors$ action (cf.\
\cite{peter5}).

\begin{definition} \label{maxmin}
We define the maximum and minimum discrepancies $M$ and $m$ for $(Y,K,\mathfrak{t})$ to be the two integers 

\begin{equation}
M = \max _{x\in \widehat{HFK}(Y,K,j)} \gr (x) - j \quad \quad \mbox{ and } \quad \quad  m = \min _{x\in \widehat{HFK}(Y,K,j)} \gr (x) - j 
\end{equation}

\end{definition}

For a knot $K\subset Y$, we let $g$ denote the integer
\[
g = max\{ j\in\zee | \widehat{HFK}(Y,K; j) \neq 0\}.
\]

\begin{theorem} \label{modulethm}
If $(M+g-2) -2k < m+ k$ then the isomorphism of theorem
\ref{morerelevantone} is a $\Lambda ^* H_1(Y)\otimes _\mathbb{Z}
\mathbb{Z}[U]$ module isomorphism.
\end{theorem}

\begin{proof} This is essentially proved in the last paragraph of the
proof of theorem 9.3 from \cite{peter1}. We sketch the argument here.
First, it follows from the proof of $\ker f_1 \cong \ker \tau_*$ (see
(\ref{SES})) that the identification in theorem \ref{morerelevantone}
is an identification of $H_1(Y)\otimes\zee[U]$-modules if the map
$f_2$ sends any element $\xi$ of $\ker f_1$ into the group
\[
X = H_*(C\{i\geq 0 \mbox{ and } j<k\}).
\]
This can be arranged under the hypotheses of the theorem by using the
degree shift formula, as follows. We impose an absolute grading on
$HF^+(Y_0,k)$ by using its identification with $\ker F$ followed by
the identification of the latter as a subgroup of $H_*(C\{i\geq 0 \mbox{ or }
j\geq k\})$ via theorem \ref{hello}. Now note that the top-degree
element of $\ker f_1$ has degree at most $M+g-2$, and compared with
$f_1$, we know $f_2$ lowers degree by at least $2k$. Furthermore, any
element $\eta\in H_*(C\{i\geq 0\})$ having degree less than $m+k$ is
necessarily contained in $X$. Hence, given $\xi\in \ker f_1$ we are
assured that $f_2(\xi)\in X$ whenever $(M+g-2)-2k< m+k$.

The remaining steps in the identification $HF^+(Y_0, k) \cong
H_*(C\{i\geq 0 \mbox{ or } j\geq k\})$ have already been seen to be
module isomorphisms.
\end{proof} 

\begin{remark} Typically the absolute grading $gr(x)$ is a nonintegral
rational number, so both $M$ and $m$ are a priori rational. However,
for the purposes of applying theorem \ref{modulethm} we may replace this
$\cue$ grading by an arbitrary $\zee$ grading (compatible with the
relative grading), since both sides of the inequality are thereby
shifted by the same amount.
\end{remark}

The 3-manifold $Y_0(K)$ carries an extra generator $\mu\in
H_1(Y_0(K);\zee)$ coming from the meridian of $K$. The above theorem does
not describe the action of this class, but we have the following:

\begin{lemma} Suppose the homomorphism $G$ in \eqref{LES} is trivial
(for example, if $HF^+_{red}(Y) = 0$ and $k \neq 0$). Then the action
of $\mu$ on $HF^+(Y_0(K), Q^{-1}([k]))$ is trivial.
\end{lemma}

\begin{proof} Since $\mu$ is homologous to a torsion class in $Y_n$
via the cobordism connecting $Y_0$ to $Y_n$, we have
\[
H(\mu\cdot \xi) = 0
\]
for $\xi\in HF^+(Y_0 , Q^{-1}([k]))$. The result follows since when $G
= 0$, $H$ is injective.
\end{proof}

\subsection{Calculational approach}

We sketch here the basis for the calculations to follow.

Let $K\subset Y$ be a nullhomologous knot in $Y$.
We will need to refer to various spin$^c$-structures on the triple of
manifolds $Y$, $Y_0$ and $Y_n$. The standard 2-handle cobordisms
connecting these manifolds, and a choice of Seifert surface $F$ for $K$,
provide identifications
\[
Spin^c(Y_0) \cong Spin^c(Y) \times \mathbb{Z} \quad \quad \quad
Spin^c(Y_n) \cong Spin^c(Y) \times \mathbb{Z}_n. 
\]
For example, $\tilde{\s}\in Spin^c(Y_0(K))$ is identified with
$(\s,k)\in Spin^c(Y)\times \zee$ if $\tilde{\s}$ is cobordant to $\s$
by a \spinc structure $\R$ on the cobordism connecting $Y$ and
$Y_0(K)$ having $\langle c_1(\R), [\hat{F}]\rangle = 2k$, where
$\hat{F}$ denotes the Seifert surface capped off by the core of the
2-handle in the cobordism. For a spin$^c$-structure $\mathfrak{s}\in
Spin^c(Y)$ we will write $\mathfrak{s} _k$ for the spin$^c$-structure
on either $Y_0$ or $Y_n$ corresponding to $(\mathfrak{s},k)$ under the
above identifications. This is consistent with the notation from the
beginning of this section.

According to the results above, in order to find $HF^+(Y_0(K), \s_k)$
we must calculate the homology of the quotient complex $C\{i<0 \mbox{
and } j\geq k\}$ of $CFK^\infty(Y, K)$ (where the filtration on the
latter group is defined using our fixed Seifert surface $F$). To do
so, we make use of a second filtration on $CFK^\infty$, namely
\[
{\mathcal F}': [\x,i,j]\mapsto i+j.
\]
Thinking of this as a filtration on $CF^\infty$, we get a spectral
sequence converging to $HF^\infty(Y,\s)$. The $E^1$ term of this
sequence may be identified with $\widehat{HFK}(Y,\s)\otimes
\zee[U,U^{-1}]$ (see lemma 3.6 in \cite{peter1}).  

Similarly, by restricting ${\mathcal F}'$ to $C\{i<0\mbox{ and }j\geq k\}$ we
obtain a spectral sequence whose $E^1$ term is
\begin{equation} \label{three1term}
Y(g,d) = \bigoplus ^d _{i= 0 } \widehat{HFK}(Y,K,g-i) \otimes
_\mathbb{Z} \mathbb{Z}[U] /U^{d+1-i}.
\end{equation}
Here we define $d = g-1-k$, and as above we let $$g=g_{\mathfrak{s}}
= \max \{ j \, | \widehat{HFK}(Y,K,\mathfrak{s},j) \ne 0 \}$$
Our method of calculation in the following sections will be to compute the
group $Y(g,d)$ and the $d_1$ differential in the spectral sequence. An
explicit calculation of the homology follows in each case, and luckily
in each case the spectral sequence collapses at the $E^2$ stage so
that no further work is necessary.

The question remains how to find the first differential $d_1$ on
$Y(g,d)$. From the structure of the filtration ${\mathcal F}'$, we see
that $d_1$ is comprised of the sum of ``horizontal'' and ``vertical''
components, the first of which decreases $i$ by one, and the second
decreases $j$ by one. To understand the vertical differential, note
that ${\mathcal F}'$ restricts to the complex $\widehat{CF} = C\{i = 0\}$
as the usual filtration induced by the knot $K$. Therefore, the
corresponding $d_1$ differential must agree with the $d_1$
differential in the spectral sequence calculating $\widehat{HF}(Y,\s)$
from $\widehat{HFK}(Y,K,\s)$, which in our examples will be a simple
matter to describe. The horizontal differential can be obtained
similarly using the fact that $C\{j = 0\}$ also is a filtered version
of $\widehat{CF}$, or alternatively by formal calculation using the
fact that $d_1^2 = 0$. More details will be given in the
following sections, as need arises.

\section{Notational shorthand and a model calculation}
\label{notationshort}

We begin our study of mapping tori by reviewing the calculation in
\cite{peter1} for the case of a single nonseparating Dehn twist. The 
connected sum theorem for $\widehat{HFK}$, which we frequently rely on in the 
remainder of the paper, can be found in section 7 of \cite{peter1}. 
With the notation used there, it states (see \cite{peter1} for more details):
$$ \widehat{HFK}(Y_1\# Y_2, K_1\#K_2, \mathfrak{t}_3) \cong 
\bigoplus _{\mathfrak{t}_1+ \mathfrak{t}_2 = \mathfrak{t}_3}
H_*\left( \widehat{CFK}(Y_1,K_1,\mathfrak{t}_1) \otimes \widehat{CFK}(Y_2,K_2,\mathfrak{t}_2)\right) $$ 

\begin{figure}[htb!] \small
\centering
\includegraphics[width=4cm]{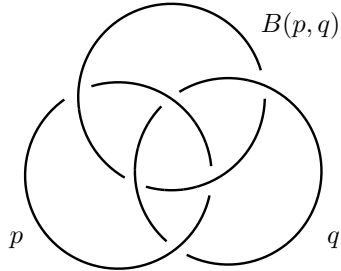}
\put(-25,90){$B(p,q)$}
\put(-120,10){$p$}
\put(0,10){$q$}
\caption{The  knot $B(p,q)$.} \label{Boromean}
\end{figure}

The first step is to realize this mapping torus $M(t_\gamma)$ as the
result of 0-surgery on some knot in a 3-manifold $Y$. Throughout, we
use $\Sigma_g$ to denote a closed surface of genus $g$. Recall that
the 3-torus $\Sigma_1\times S^1$ can be realized as the result of
0-surgery on each component of the Borromean rings. Alternatively, if
we perform 0-surgery on two components, the remaining component
describes a knot in the connected sum $S^1\times S^2\#S^1\times S^2$.
Write $B(0,0)$ for this knot: then $\Sigma_1\times S^1$ is given by
0-surgery on $B(0,0)$. To produce the mapping torus of a right-handed Dehn
twist about a nonseparating curve $\gamma$ in $\Sigma_1$, we may
modify the surgery picture by adding a $-1$-framed meridian to one of
the surgery circles in the picture for $2(S^1\times S^2)$. Blowing
down this circle gives a knot $B(1,0)$ in $S^1\times S^2$, and
performing 0-surgery on $B(1,0)$ produces the desired mapping torus 
(see figure \ref{Boromean}).

For general $g$, we have a surgery description of $\Sigma_g\times S^1$
as 0-surgery on the connected sum of $g$ copies of $B(0,0)$, which is
a knot in the connected sum $2g(S^1\times S^2)$. Just as above, the
mapping torus $M(t_\gamma)$ of a Dehn twist on $\gamma$ is obtained by
0-surgery on $B(1,0)\#(g-1)B(0,0)$. 

It is proved in \cite{peter1} that the knot Floer homologies of
$B(0,0)$, $B(1,0)$, and $B(-1,0)$ (the latter corresponding to a
left-handed Dehn twist) are trivial unless we use the torsion $\spinc$
structure on the underlying 3-manifold, and in this case we have:
\begin{eqnarray}
\widehat{HFK}(2(S^1\times S^2), B(0,0) , j) &\cong&
\Lambda^{j+1}H^1(\Sigma_1)_{(j)}\nonumber\\
\widehat{HFK}(S^1\times S^2, B(1,0), j) &\cong&
\Lambda^{j+1}H^1(\Sigma_1)_{(j-\frac{1}{2})}\label{HFKB}\\
\widehat{HFK}(S^1\times S^2, B(-1,0), j) &\cong&
\Lambda^{j+1}H^1(\Sigma_1)_{(j+\frac{1}{2})}.\nonumber
\end{eqnarray}
The spectral sequence for $\widehat{HF}$ coming from $B(0,0)$ has no
nontrivial differentials for dimensional reasons, while those coming from
$B(1,0)$ and $B(-1,0)$ have nontrivial differentials only at the $E^1$
stage. In any case, it follows from the connected sum theorem for
$\widehat{HFK}$ 
that as relatively graded groups all three of the
knots $gB(0,0)$, $B(1,0)\#(g-1)B(0,0)$ and $B(-1,0)\#(g-1)B(0,0)$ have
knot Floer homology $\widehat{HFK}$ isomorphic (as a relatively graded
group) to $\Lambda^*H^1(\Sigma_g)$, where the filtration on
$\Lambda^kH^1(\Sigma_g)$ is given by $k-g$.

The $E^1$ term of the spectral sequence for $C\{i<0\mbox{ and }j\geq
k\}$ is given in general by (\ref{three1term}), which in this case
takes the form
\begin{equation} \label{Xgddef}
X(g,d) = \bigoplus ^d _{i = 0} \Lambda ^{2g-i} H^1(\Sigma _g) \otimes
_\mathbb{Z} \mathbb{Z}[U] /U^{d+1-i}.
\end{equation}
We will continue to use the notation $X(g,d)$ for the above group in
all of what follows. Note that according to \cite{macdonald}, $X(g,d)$
is isomorphic to  $H_*(Sym^d(\Sigma _g))$, the homology of the $d$-fold
symmetric product of $\Sigma _g$.

It is convenient to picture $X(g,d)$ as an array of groups:
\begin{equation}
\begin{array}{cccccc}
 & & & & & \Lambda^{2g}H^1 \\
 & & & &\Lambda^{2g}H^1\otimes U & \Lambda^{2g-1}H^1\\ &&& &&\\
 & &\cdots & \cdots & \cdots & \vdots\\
 & &&&&\\
&\Lambda^{2g}H^1\otimes U^{d} & \Lambda^{2g-1}H^1\otimes U^{d-1} &
\cdots & \cdots & \Lambda^{2g-d}H^1\end{array}
\label{arraypic}
\end{equation}
\begin{definition} \label{gradofsigma}
We define a grading on $X(g,d)$ by the formula
$\deg(\omega\otimes U^j) = \deg (\omega) - g - 2j$, where $\deg(\omega)$
indicates the degree of the homogeneous element $\omega$ in the exterior
algebra. We will use the notation $X(g,d){[n]}$ to indicate the graded
group $X(g,d)$ with the grading shifted up by $n$, and similar
notation for other groups with obvious gradings.
\end{definition}

This definition is consistent with the convention from the
introduction; the grading here is the one obtained from the absolute
grading on $CFK^\infty(gB(0,0))$.
Note that this grading will not carry over to a
grading on $HF^+$ for the mapping tori we are considering (the latter
groups do not carry integral gradings), but it will correspond to a
relative cyclic grading in those groups.

There is an action of $H_1(\Sigma_g)$ on $X(g,d)$ given as follows.
Let $\gamma\in H_1(\Sigma_g)$, and define $D_{\gamma}$ by 
\begin{equation} \label{Dgamma}
D_\gamma(\omega\otimes U^j) = \iota_\gamma\omega\otimes U^j +
PD(\gamma)\wedge\omega\otimes U^{j+1}.
\end{equation}
In terms of the array of groups \eqref{arraypic}, $D_\gamma$ maps an element in
one of the groups in the array into the two groups immediately to the
left and immediately below the starting point. An easy check reveals
that $D_\gamma D_\gamma =0$, thus for a fixed curve $\gamma$,
$D_\gamma$ defines a differential on $X(g,d)$.

For $g\ge 2$ we write $\Sigma _g = \Sigma_1 \# \Sigma _{g-1}$.
This induces a splitting
$
\Lambda^* H^1(\Sigma_g) = \Lambda^*_+H^1(\Sigma_g) \oplus
\Lambda^*_-H^1(\Sigma_g) 
$
where
\begin{eqnarray*}
\Lambda^*_+H^1(\Sigma_g) &=& (\Lambda^0H^1(\Sigma_1) \oplus \Lambda^2H^1(T)
)\otimes \Lambda^*H^1(\Sigma_{g-1})\\
\Lambda^*_-H^1(\Sigma_g) &=& \Lambda^1H^1(\Sigma_1)\otimes
\Lambda^*H^1(\Sigma_{g-1}).
\end{eqnarray*}
This in turn induces a splitting $X(g,d) = X_+(g,d) \oplus X_-(g,d)$.
Let $\gamma \subset \Sigma_1$ denote a fixed simple closed curve on
$\Sigma_1$. Using the differential $D_\gamma$ from \eqref{Dgamma}, we
define two new differentials $D^\pm _\gamma$ as
\begin{equation} \label{Dgammapm}
D^\pm _\gamma | _{X_\pm (g,d)} = 0\quad \quad \quad \quad
\quad D^\pm _\gamma | _{X_\mp (g,d)} = D_\gamma 
\end{equation}

\begin{remark} Our notation here differs slightly from \cite{peter1},
where the roles of $D_\gamma^+$ and $D_\gamma^-$ are exchanged. The
conventions here ensure that ``$+$'' is used for a right-handed twist
throughout.
\end{remark}

\begin{lemma} For $K_\pm$ the knot $B(\pm 1, 0)\#(g-1)B(0,0)$, the
$d_1$ differential on $X(g,d)$ in the spectral sequence for
$HF^+(M(t_\gamma^{\pm}), k)$ is given by $D_\gamma^\pm$. All subsequent
differentials are trivial.
\end{lemma}

The claim about further differentials vanishing follows immediately
from dimensional considerations. This lemma is proved in
\cite{peter1}; the fact that $d_1 = D_\gamma^\pm$ is a straightforward
check using the description of $d_1$ at the end of the
previous section. In particular, we can arrange the isomorphism
\[
\widehat{HFK}((2g-1)S^1\times S^2, B(\pm 1, 0)\#(g-1)B(0,0)) \cong
\Lambda^*H^1(\Sigma_g)
\]
in such a way that the $d_1$ differential in the spectral sequence for
$\widehat{HF}$ is trivial on $\Lambda^*_+$ (resp $\Lambda^*_-$) and given by
contraction with $\gamma$ on $\Lambda^*_-$ (resp $\Lambda^*_+$). (This in
turn is clear from the expressions (\ref{HFKB}) and the connect-sum
theorem for $\widehat{HFK}$.) This shows that the vertical
differential is the vertical component of $D_\gamma^\pm$; to check the
horizontal differential, one can see that it must have the claimed
form by enforcing the equation $d_1^2 = 0$.

Our object now is to compute the homology of $(X(g,d), D_\gamma^\pm)$.
We focus on the case of a right-handed twist, meaning that the
relevant differential is $D_\gamma^+$. 

First we obtain an alternate
description of $X_+(g,d)$. Let us write $c^*\in H^1(\Sigma_1)$ for the
Poincar\'e dual of $\gamma$ and $c\in H^1(\Sigma_1)$ for the Poincar\'e
dual of a geometric dual generator in homology. In other words, we
have
\begin{equation}
\begin{array}{ll}
\iota_\gamma c = 1 \qquad & PD(\gamma)\wedge c = c^*\wedge c \\
\iota_\gamma c^* = 0 & PD(\gamma)\wedge c^* = 0.
\end{array}
\label{cconv}
\end{equation}
Thus an element $\omega\otimes U^j\in X_+(g,d)$ ($\omega$ a
monomial) has the property that either $\omega\in
\Lambda^*H^1(\Sigma_{g-1})$ or $\omega$ is of the form $c\wedge c^*
\wedge\lambda$ for $\lambda\in \Lambda^*H^1(\Sigma_{g-1})$. Therefore:
\begin{eqnarray}
X_+(g,d) &=& \left(\bigoplus_{i = 0}^d \Lambda^{2g-i}H^1(\Sigma_{g-1})
\otimes \zee[U]/U^{d-i+1} \right)\nonumber\\ &&\hspace*{2em} \oplus \left(
c\wedge c^* \wedge \bigoplus_{i = 0}^d
\Lambda^{2g-i-2}H^1(\Sigma_{g-1})\otimes \zee[U]/U^{d-i+1}\right)
\nonumber\\
&=& \left( \bigoplus_{j = 0}^{d-2} \Lambda^{2g-2-j}H^1(\Sigma_{g-1})
\otimes \zee[U]/U^{(d-2) - j + 1}\right) \nonumber\\ &&
\hspace*{2em}\oplus \left( c\wedge c^* \wedge \bigoplus_{i = 0}^d
\Lambda^{2g-i-2}H^1(\Sigma_{g-1})\otimes
\zee[U]/U^{d-i+1}\right)\nonumber\\
&\cong& X(g-1,d-2)[-1] \oplus c\wedge c^* \wedge X(g-1,d)[-1].\label{X+form}
\end{eqnarray}
(Note that in the first summand of the first line, the terms with $i =
0$ and $i = 1$ vanish, so the second line follows by replacing $i$ by
$j = i + 2$.)

Now we consider the homology of $(X(g,d), D_\gamma^+)$. First, note
that $X_+(g,d)$ consists of cycles by definition. To see the homology
coming from $X_+(g,d)$, note that if $\omega\in
\Lambda^*H^1(\Sigma_{g-1})$ then from (\ref{cconv}):
\begin{eqnarray}
D_\gamma^+(c\wedge \omega\otimes U^j) &=& \omega \otimes
U^j - c\wedge c^* \wedge\omega\otimes U^{j+1}\label{Dc}\\
D_\gamma^+(c^*\wedge\omega\otimes U^j) &=& 0\label{Dc*}
\end{eqnarray}
From (\ref{Dc}), it follows that the portion of the homology generated
by $X_+(g,d)$ is simply $X_+(g,d)$ with the relation that
$\omega\otimes U^j - c\wedge c^*\wedge\omega\otimes U^{j+1}=0$ for
$\omega\in\Lambda^*H^1(\Sigma_{g-1})$ and $j\geq 0$. In terms of the
decomposition (\ref{X+form}) of $X_+(g,d)$, this means that the group
$X(g-1,d-2)$ appearing in (\ref{X+form}) is identified with its image
in $c\wedge c^*\wedge X(g-1,d)$ under the map 
$f( \omega\otimes U^j) = c\wedge c^*\wedge\omega\otimes U^{j+1}$. 
Since the latter map
is injective, the homology will contain a copy of $X(g-1,d-2)$. The
remainder of the homology coming from $X_+(g,d)$ is given by the
cokernel of $f$ in $c\wedge c^*\wedge X(g-1,d)$, modulo any
boundaries. Now we can think of $f$ as mapping one step diagonally up
and to the left in the array (\ref{arraypic}) (with $g$ replaced by
$g-1$); the cokernel of $f$ is then the ``bottom row'' together with
the ``rightmost column.'' An element of the rightmost column is never
a boundary, while an element $c\wedge c^*\wedge\omega\otimes U^j \in
X_+(g,d)$ that lives on the bottom row (but not in the rightmost
column) is the image of $-c\wedge\omega\otimes U^{j-1}$ under $D_\gamma^+$.
Thus the cokernel of $f$ contributes the group
\[
c\wedge c^* \wedge \bigoplus_{i = 0}^d \Lambda^{2g-2-i}
H^1(\Sigma_{g-1}) \cong \bigoplus_{i = 0}^d
\Lambda^{2g-2-i}H^1(\Sigma_{g-1})_{(g-i)}
\]
to the homology. In fact we can rewrite the contribution of $X_+(g,d)$ in
another way, as follows.

So far, $X_+(g,d)$ contributes 
\[
X(g-1,d-2)[-1] \oplus \bigoplus_{i = 0}^d
\Lambda^{2g-2-i}H^1(\Sigma_{g-1})_{(g-i)}.
\]
Comparing this with the array (\ref{arraypic}), we see that 
\[
X(g-1,d-2)[-1] \oplus\bigoplus_{i = 0}^{d-1}
\Lambda^{2g-2-i}H^1(\Sigma_{g-1})_{(g-i)} \cong X(g-1, d-1)[1].
\]
Hence the contribution from $X_+(g,d)$
can be written as
\[
X(g-1,d-1)[1] \oplus \Lambda^{2g-2-d}H^1(\Sigma_{g-1})_{(g-d)}.
\]

Consider now the contribution from $X_-(g,d)$. According to (\ref{Dc}) and
(\ref{Dc*}), and the fact that no element of $X_-(g,d)$ can be a
boundary, the homology coming from $X_-(g,d)$ is isomorphic to the
group generated by elements of the form $c^*\wedge\omega\otimes U^j$,
$\omega\in \Lambda^*H^1(\Sigma_{g-1})$. This is
\begin{eqnarray*}
\langle c^*\wedge\omega\otimes U^j \rangle &=& c^*\wedge\bigoplus_{i
=0}^{d-1} \Lambda^{2g-2 - i}H^1(\Sigma_{g-1}) \otimes \zee[U]/U^{d-i} \\
&=& c^*\wedge X(g-1,d-1)[-1] \cong X(g-1,d-1).
\end{eqnarray*}
This completes the calculation of $H_*(X(g,d), D^+_\gamma)$; the case
of $H_*(X(g,d), D^-_\gamma)$ is treated similarly. To summarize, we
have proved the case of the following theorem relevant to a
right-handed Dehn twist $t_\gamma$, which proves the case $n = 1$ of
theorem \ref{mainmain}. (Note that when $n = 1$ there is a unique
\spinc structure $\s_k\in{\mathcal S}_k$.)
\begin{theorem} \label{homologyofxgd}
Let $\Sigma _g = \Sigma_1 \# \Sigma _{g-1}$ be a genus $g$ surface with
$g\ge 2$ and let $\gamma \subset \Sigma_1$ be a simple closed curve. With
$X(g,d)$ and $D^\pm _\gamma$ as defined by \eqref{three1term} and
\eqref{Dgammapm}, the homology groups $H_*(X(g,d), D^\pm _\gamma)$ are
isomorphic as relatively graded groups to $HF^+(M(t_\gamma^{\pm 1}),
\s_k)$, where $d = g-1-|k|$. Explicitly, we have:
$$
\begin{array}{l} 
H_*(X(g,d), D^+ _\gamma) =  \cr 
\quad \quad \quad \quad X(g-1,d-1)\oplus X(g-1,d-1){[1]} \oplus \Lambda ^{2g-2-d}H^1 (\Sigma _{g-1})_{(g-d)}  \cr
\cr 
H_*(X(g,d), D^- _\gamma) =  \cr 
\quad \quad \quad \quad X(g-1,d-1)\oplus X(g-1,d-1){[-1]} \oplus \Lambda ^{2g-2-d}H^1 (\Sigma _{g-1})_{(g-d)} 
\end{array}
$$
with $X(g-1,\cdot){[\ell \, ]}$ as in definition \ref{gradofsigma}.
\end{theorem}

The calculation for a left-handed twist is dual to the one just performed. While the calculation preceding 
theorem \ref{homologyofxgd} focused on the case $k>0$, the case of $-k$ yields the same result. This 
can be seen from formula \eqref{negativeconversion} together with the observation that the \lq\lq background\rq\rq
torsion \spin $\mathfrak{s}_0$ on $Y=\#(2g-1)(S^1\times S^2)$ is its own conjugate.

\section{Multiple Dehn twists along a non-separating curve} \label{five5}
In this section we calculate $HF^+$ for the mapping torus $M(t_\gamma^n)$
where $t_\gamma$ denotes the right-handed Dehn twist along a
nonseparating curve $\gamma \subset \Sigma _g$ as usual. As in the previous
section, our strategy will be to exhibit $M(t_\gamma^n)$ as the result of
zero surgery along a knot and then to use the tools developed in
section \ref{prelim}.

The mapping torus $M(t_\gamma^n)$ is obtained as zero framed surgery along
the knot
\begin{equation}\label{ntwistKdef}
K = B(n,0) \# (g-1) B(0,0) 
\end{equation}
in $ Y = L(n,1)\#(2g-1) (S^1\times S^2)$. Here, as in \cite{peter1},
if $p,q\in \cue\cup\{\infty\}$ then we can perform surgery on two
components of the Borromean rings with coefficients $p$ and $q$: the
third component is the knot $B(p,q)$ (see figure \ref{Boromean}).

We first calculate $\widehat{HFK}(Y, K)$. To do so, note that
$H^2(Y;\zee) \cong \zee_n\oplus\zee^{2g-1}$, and in fact we can write a
\spinc structure on $Y$ as a connected-sum $\T_i\#\s$ where $\s\in
Spin^c((2g-1)(S^1\times S^2))$ and $\T_i\in Spin^c(L(n,1))$, $i\in\zee_n$ is
defined using the conventions from section \ref{prelim} (ie, using
the standard 2-handle cobordism from $S^3$ to $L(n,1)$). We write $\s_0$
for the torsion \spinc structure on $\#(2g-1)S^1\times S^2$; we also
assume that $n>0$ since the case $n<0$ is treated in much the same
way. First we give the calculation in the case $g = 1$.
\begin{lemma} \label{multiknot1}
The knot Floer homology $\widehat{HFK}(L(n,1)\#(S^1\times S^2), B(n,0)
, \mathfrak{t}_i\# \mathfrak{s} )$ is zero unless $\s = \s_0$. In this
case, if $i\ne 0$ we have isomorphisms of relatively graded groups
\begin{align}\label{boringmulti}
\widehat{HFK}(L(n,1)\#(S^1\times S^2),
B(n,0), \mathfrak{t}_i\# \mathfrak{s}_0 , j ) \cong \left\{
\begin{array}{cl} H^*(S^1;\zee) \quad & j = 0 \cr 0 & \mbox{otherwise}
\end{array} 
\right. 
\end{align}
while for $i=0$ we have
\begin{align}
\widehat{HFK}(L(n,1)\#(S^1\times S^2), B(n,0),  
\mathfrak{t}_0\# \mathfrak{s}_0 , j ) \cong
\Lambda^{j+1}H^1(\Sigma_1;\zee),
\end{align}
where $\Sigma_1$ denotes a genus-1 surface. The spectral sequence for
$\widehat{CF}$ induced by $K$ collapses at the $E^1$ level for $i\ne
0$ and collapses at the $E^2$ level for $i=0$. In the latter case, the
$d_1$ differential is trivial except for its restriction to the $j=0$
summand
\[
d_1 : \Lambda^1H^1(\Sigma_1)\to \Lambda^0H^1(\Sigma_1)
\]
which is surjective.

\end{lemma}

\begin{proof}
The starting point for this calculation is the surgery long exact
sequence in knot Floer homology (see \cite{peter3} and \cite{peter1}):
\begin{align} \label{LES3}
\cdots\rightarrow \widehat{HFK}(S^1 \times S^2 ,  & B(\infty,0),\s,j) 
\rightarrow \widehat{HFK}(\#^2(S^1\times S^2), B(0,0),
[\s_i]\#\s,j ) \rightarrow \cr 
\rightarrow &\widehat{HFK}( L(n,1)\#
(S^1\times S^2) , B(n,0), \mathfrak{t}_i\#\s,j ) \rightarrow \cdots
\end{align}
Of course $B(\infty , 0)$ is just the unknot in $S^1\times S^2$. When
$\s \neq \s_0$, both the first and second term in (\ref{LES3}) vanish
by the connected sum theorem for $\widehat{HFK}$, so the first claim
of the lemma is immediate. For $\s = \s_0$, there are several cases to
consider, depending on $\mathfrak{t}_i$ as well as the filtration
index $j$.

When $i\ne 0$ then the second term in \eqref{LES3} is zero and the
first term is only nonzero for $j=0$. Indeed, since the unknot induces
the trivial filtration on $\widehat{CF}$, the knot Floer homology of
the unknot in $Y$ is isomorphic to $\widehat{HF}(Y)$ supported in
filtration level zero. Since $\widehat{HF}(S^1\times S^2; \s_0)\cong
H^*(S^1;\zee)$ as relatively graded groups, (\ref{boringmulti}) follows.

If $i=0$ and $j\ne 0$, then the first term in \eqref{LES3} is zero
and the second term is nonzero only if $j=\pm 1$ (cf.\ equation
\eqref{HFKB}). Since $\widehat{HFK}(B(0,0),j)\cong \zee$ when $j = \pm
1$, the statement follows in this case.

In the remaining case $i=j=0$, (\ref{LES3}) easily shows that
$\widehat{HFK}(B(n,0))$ can be nonzero only in two adjacent degrees:
let us call them $a$ and $b$. Then the sequence becomes

\begin{align} \nonumber
0  \rightarrow \widehat{HFK}_a( L(n,& 1)\#(S^1 \times S^2), B(n,0) , \mathfrak{t}_0\#\mathfrak{s}_0 ,0 )    \rightarrow \mathbb{Z}_{\left( \frac{1}{2}\right)} \rightarrow 
\mathbb{Z}^2_{(0)} \rightarrow  \cr
\rightarrow & \widehat{HFK}_b(L(n,1)\#(S^1\times S^2), B(n,0),
\mathfrak{t}_0\#\mathfrak{s}_0,0 ) \rightarrow \mathbb{Z}_{\left( -
\frac{1}{2}\right)} \rightarrow 0
\end{align}
and it is an exercise to see that the map $\mathbb{Z}_{\left(
\frac{1}{2}\right)} \rightarrow \mathbb{Z}^2_{(0)}$ is injective. This
completes the proof, as the statements about the spectral sequence follow
easily from dimensional considerations.
\end{proof}

For the case of general $g$ in \eqref{ntwistKdef}, the connected sum
theorem for $\widehat{HFK}$ combined with the previous lemma yields:
\begin{lemma}  \label{knotfloer3}
The Heegaard knot Floer homology $\widehat{HFK}( Y, K,\mathfrak{t})$ is zero
unless $\mathfrak{t} = \mathfrak{t}_i\#^{2g-1} \mathfrak{s}_0$ for some
$i$ and in this case is given by
\begin{equation} \label{explicit1}
\widehat{HFK}(Y,K, \mathfrak{t}_i\#^{2g-1} \mathfrak{s}_0  )  \cong \left\{ 
\begin{array}{ll}
\Lambda ^*(H^1(S^1;\mathbb{Z}) \oplus   H^1(\Sigma _{g-1} ;\mathbb{Z})) \quad \quad &  i \ne 0 \cr
\Lambda ^*H^1(\Sigma_{g} ;\mathbb{Z}) & i = 0 \end{array} \right.
\end{equation}
The spectral sequence for $\widehat{CF}$ induced by $\widehat{CFK}$
collapses at the $E^1$ level for $i\ne 0$ and it collapses at the
$E^2$ level for $i=0$. The $d_1$ differential in the latter case is
trivial except on $H^1(\Sigma_1)$ where it is contraction with
$\gamma$ (a generator of $H^1(\Sigma_1)$).
\end{lemma} 
Let $Y(g,d)$ be as in \eqref{three1term}. In the case $i=0$ we see from
\eqref{explicit1} that in fact $Y(g,d) = X(g,d)$. Thus the splitting
$\Sigma _g = \Sigma_1\# \Sigma _{g-1}$ induces the decomposition
$Y(g,d) = X(g,d) = X_+(g,d) \oplus X_-(g,d)$ as in section
\ref{notationshort}; the $d_1$ differential from lemma \ref{knotfloer3}
gives rise to a differential $D$ on $Y(g,d)$ which can easily be seen
to be given by $D^+_\gamma$. Just as in the case $n = 1$ there can be no
differentials beyond $d_1$ in the spectral sequence for $H_*(C\{i<0
\mbox{ and } j\geq k\})$, so according to theorem
\ref{morerelevantone} the Floer homology $HF^+(M(t_\gamma^n))$ is
given as the homology $H_*(Y(g,d),D)$ (with $D = 0$ if $i\ne 0$). The
assumptions $HF^+_{red}(Y) = HF^+_{red}(Y_n)=0$ of theorem
\ref{morerelevantone} are readily verified. For the case $i\neq 0$,
lemma \ref{knotfloer3} shows that there can be no $d_1$ differential
in the spectral sequence for $\widehat{HF}$---indeed since the right
hand side of \eqref{explicit1} is equal to $\widehat{HF}(Y,\T_i)$ the
spectral sequence collapses at the outset.

For stating the main theorem of this section, we will use the notation
$\mathfrak{t}_{i,k}$ to denote the spin$^c$-structure
$(\mathfrak{t}_i\#\mathfrak{s}_0 ,k) \in Spin^c(Y)\times
\mathbb{Z} \cong Spin^c(M(t_\gamma^n))$.
\begin{theorem} \label{main2} 
Let $M(t_\gamma^n)$ be the mapping torus of $n$ right-handed Dehn twists
along a non-separating curve $\gamma \subset \Sigma _g$, $g\ge 2$,
where $n\in\zee\setminus\{0\}$. Pick an integer $k$ with $0<|k|\leq g-1$.
Then the Heegaard Floer homology of $M(t_\gamma^n)$ for the
spin$^c$-structure $\mathfrak{t}_{i,k}$ is given as a relatively
graded group by
$$ HF^+(M(t_\gamma^n) , \mathfrak{t}_{i,k} ) \cong \left\{
\begin{array}{ll}
H^*(S^1) \otimes X(g-1,d-1)   \quad & i\ne 0  \cr
& \cr 
 H_*(X(g,d), D_\gamma ^{s (n)} ) & i =0
\end{array}
\right.  $$
where $s (n) = \mbox{sign}\, (n)$ and $d = g-1-|k|$. If in addition
$g<3k+2$, the above isomorphisms are $\Lambda ^* H_1(Y)\otimes
_\mathbb{Z} \mathbb{Z}[U]$ module isomorphisms.
\end{theorem}
The claims about the module structure follow easily from
lemma \ref{knotfloer3} and theorem \ref{modulethm}. The groups $H_*(X(g,d),
D_\gamma ^\pm)$ have explicitly been calculated in theorem
\ref{homologyofxgd} and are given by
$$\begin{array}{l} 
H_*(X(g,d),D^+ _\gamma) =  \cr 
\quad \quad \quad  X(g-1,d-1)\oplus X(g-1,d-1){[1]} \oplus \Lambda ^{2g-2-d}H^1 (\Sigma _{g-1})_{(g-d)}  \cr
 \cr 
H_*(X(g,d),D^- _\gamma) = \cr 
\quad \quad \quad  X(g-1,d-1)\oplus X(g-1,d-1){[-1]} \oplus \Lambda ^{2g-2-d}H^1 (\Sigma _{g-1})_{(g-d)} 
\end{array} $$

This proves theorem \ref{mainmain} for $k>0$. The case of $k<0$ follows from the observation that 
$\mathfrak{t}_{i,-k} = \overline{\mathfrak{t}_{n-i,k}}$. In particular, 
$\mathfrak{t}_{0,-k} = \overline{\mathfrak{t}_{0,k}}$. Since 
$HF^+(M(t_\gamma^n) , \mathfrak{t}_{i,k} )$ (in the case when $k>0$) only depends on if $i=0$ or 
$i\ne 0$ (a property which is preserved under conjugation of the \spin ), the theorem for $k<0$ now 
readily follows from formula \eqref{negativeconversion}. 


\section{Multiple Dehn twists along a transverse pair of curves} \label{six6}

Let
$\gamma,\delta \subset \Sigma_1$ be a pair of curves representing
generators of $H_1(\Sigma_1)$ with $\gamma \cap \delta$ being a single
point. Write $\Sigma _g = \Sigma_1\# \Sigma _{g-1}$ (with $g\ge 2$) and
think of $\gamma$ and $\delta$ as curves in $\Sigma_g$. We consider
the mapping torus $M(t_\gamma^mt_\delta^n)$.

The manifold $M(t_\gamma^mt_\delta^n)$ can be described as the result of
zero surgery along the knot
\begin{equation} \label{knotdef3}
 K = B(m,n)\, \#(g-1)B(0,0)
\end{equation}
inside $Y=L(m,1)\#L(n,1)\#(2g-2)(S^1 \times S^2)$. Using the standard
2-handle cobordism from $(2g-2)(S^1\times S^2)$ to
$Y$ (having two handles), we have an identification as in section
\ref{prelim}:
\[
Spin^c(Y) \cong \zee_m\oplus \zee_n\oplus Spin^c((2g-2)(S^1\times S^2)).
\]
We will use $\mathfrak{t}_{a,b}$ to denote the spin$^c$-structure
corresponding to $(a,b,\s_0)$ under this identification. To calculate the
knot Floer homology of $K$, we first calculate the knot Floer homology
of $B(m,n)$ and then use the connected sum theorem for
$\widehat{HFK}$. The discussion proceeds by considering the two cases
$m\cdot n >0$ and $m \cdot n <0$ separately.

\subsection{The case $m\cdot n >0$} \label{hellohello}

We assume that $m$ and $n$ are both positive for the moment. We find the
knot Floer homology of $B(m,n)$ from the surgery long exact sequence:

\begin{align} \label{LES2}
\cdots\rightarrow \widehat{HFK}(L(m,1)\#  S^3, B(m,\infty), \mathfrak{t}_a  )   \rightarrow  & \cr
\rightarrow \widehat{HFK}(L(m,1) \# (S^1 \times  S^2) , B(m & ,0) , [\mathfrak{t}_{a,b}]) \rightarrow \cr
\rightarrow  \widehat{HFK}(L(m,1)\# L & (n,1), B(m,n), \mathfrak{t}_{a,b} )  \rightarrow \cdots
\end{align}
In the above, $B(m,\infty)$ is the unknot in $L(m,1)$ and
$\widehat{HFK}(L(m,1)\# (S^1\times S^2) , B(m,0), \mathfrak{t}_{a,b})$
was calculated in lemma \ref{multiknot1}. From this one can see that
the sequence \eqref{LES2} essentially takes two possible forms
depending on if $b=0$ or if $b\ne 0$. In the case when $b=0$ there is
a further distinction between the cases $a=0$ and $a\ne 0$. What
follows is a case by case analysis of the sequence \eqref{LES2}.

Recall (see \cite{peter2}) that $\widehat{HF}(L(m,1),\T_a)\cong \zee$ for
all $a\in\{0,\ldots,m-1\}$, and is supported in degree $\ell_{m,a}$ where
\begin{equation} \label{ella} 
\ell_{m,a} = \frac{1}{4} \left( \frac{(2a-m)^2}{m} - 1 \right)
\end{equation}

\subsubsection{The case $b\ne 0$.}

In this case the middle term of \eqref{LES2} is zero while the first
term is nonzero only when $j=0$:
$$ 0 \rightarrow  \widehat{HFK}(L(m,1)\# L(n,1), B(m,n), \mathfrak{t}_{a,b}, 0 )  \rightarrow   
\mathbb{Z}_{(\ell_{m,a})}   \rightarrow  0  $$

\subsubsection{The case $b=0$ and $a\ne0$.}
The long exact sequence \eqref{LES2} in this case becomes (for $j=0$)
\begin{align} \nonumber
& 0 \rightarrow  \mathbb{Z}_{\left( \ell_{m,a} + \frac{1}{2} \right) }
\rightarrow \widehat{HFK}_{(\ell_{m,a} + \ell_{n,0})} (L(m, 1)\#
L(n,1), B(m,n), \mathfrak{t}_{a,0},0) \rightarrow \cr
& \mathbb{Z}_{(\ell_{m,a})} {\rightarrow}\, \mathbb{Z}_{\left( \ell_{m,a} - \frac{1}{2}
\right) } {\rightarrow}\, \widehat{HFK}_{(\ell_{m,a} +
\ell_{n,0} - 1 )} (L(m,1)\# L(n,1), B(m,n), \mathfrak{t}_{a,0},0) 
{\rightarrow}\, 0
\end{align}
It is easy to check that $\mathbb{Z}_{( \ell_{m,a})} \rightarrow 
\mathbb{Z}_{\left( \ell_{m,a} - \frac{1}{2} \right) } $ is an isomorphism. 

\subsubsection{The case $b=0$ and $a=0$.}
In this case the sequence \eqref{LES2}  splits into three exact sequences corresponding to $j=-1, 0 , 1$. In the cases $j=\pm 1$ the first term in 
\eqref{LES2} is zero which reduces the sequence to   
\begin{align} \nonumber
 0 \rightarrow \widehat{HFK}(L(m,1) \# (S^1 \times S^2), B(n,0) ,&
 [\mathfrak{t}_{0,0}],j) \rightarrow \cr \rightarrow
 \widehat{HFK}(L(m,1)\# & L(n,1), B(m,n), \mathfrak{t}_{0,0},j)
 \rightarrow 0
\end{align}
On the other hand, when $j=0$ the sequence \eqref{LES2} becomes 
\begin{align} \nonumber
0 & \rightarrow \widehat{HFK}_{(\ell_{m,0} + \ell_{n,0})} (L(m, 1)\#
L(n,1), B(m,n), \mathfrak{t}_{0,0},0) \rightarrow \mathbb{Z}_{(
\ell_{m,0})} \rightarrow \cr 
& \rightarrow  \mathbb{Z}^2_{\left( \ell_{m,0} - \frac{1}{2}
\right) } 
 \rightarrow \widehat{HFK}_{(\ell_{m,0} +
\ell_{n,0} - 1 )} (L(m,1)\# L(n,1), B(m,n), \mathfrak{t}_{0,0},0)
\rightarrow 0,
\end{align}
and one checks that $\zee_{(\ell_{m,0})}\to \zee^2_{(\ell_{m,0} -
\frac{1}{2})}$ is injective.
We summarize the above analysis in the following lemma:
\begin{lemma} \label{bmn1}
Suppose $a$ and $b$ are not both zero. Then the knot Floer homology of
$B(m,n)$ in $L(m,1)\# L(n,1)$ with $m,n>0$ in the spin$^c$-structure
$\mathfrak{t}_{a,b}$ is given by
\[
\widehat{HFK}(L(m,1)\# L(n,1), B(m,n), \mathfrak{t}_{a,b}) \cong\zee,
\]
supported in filtration level zero and degree $\ell_{m,a} + \ell_{n,b}$. 

If $a = b = 0$, then
\begin{equation} \nonumber
\widehat{HFK}(L(m,1)\# L(n,1), B(m,n), \mathfrak{t}_{0,0},j ) \cong
\left\{
\begin{array}{ll} 
\mathbb{Z}_{( \ell_{m,0} + \ell_{n,0} )} \quad
& j=1 \cr \mathbb{Z}_{( \ell_{m,0} + \ell_{n,0} -1 )} \quad
& j=0 \cr \mathbb{Z}_{( \ell_{m,0} + \ell_{n,0}-2 )} \quad &
j=-1 \cr 0 & \mbox{otherwise,} \end{array} \right.
\end{equation}
where $\ell_{m,a}$ and $\ell_{n,b}$ are given by \eqref{ella}. The
spectral sequence on $\widehat{CF}$ induced by $\widehat{CFK}$
collapses at the $E^1$ level if either $a\ne 0$ or $b\ne 0$ and it
collapses at the $E^2$ level if $a=b=0$. In the latter case, the $d_1$
differential is the surjective map from $\mathbb{Z}_{( \ell_{m,0} +
\ell_{n,0} -1 )} \rightarrow \mathbb{Z}_{( \ell_{m,0} + \ell_{n,0} -2
)} $.
\end{lemma} 
At the beginning of this section we assumed $m,n>0$. The knot Floer
groups of $B(m,n)$ for $m,n<0$ can be found from the symmetry relation
$$\widehat{HFK}_d(B(m,n), \mathfrak{t}_{a,b},j) =
\widehat{HFK}_{-d}(B(-m,-n), \mathfrak{t}_{a,b},-j)
$$
Our calculation for the knot Floer homology of $B(m,n)$ immediately
yields the knot Floer homology of the knot $K$ from \eqref{knotdef3}.
To express the result in a notation-friendly way, we'll denote by
$\Sigma _{1/2}$ the (fictitious) surface of genus $1/2$ whose
cohomology is given by $H^k(\Sigma _{1/2};\mathbb{Z}) = \mathbb{Z}$
for $ k=0,1,2$.

\begin{lemma} \label{bmn2}
For $m,n>0$ the knot Floer homology of $K$ with respect to the
spin$^c$-structure $\mathfrak{t}_{a,b}$ is
\begin{align} \nonumber
\widehat{HFK}& (Y,K,\mathfrak{t}_{a,b}) = \left\{ 
\begin{array}{ll}
\Lambda ^*H^1(\Sigma _{g-1}) \quad \quad & a\ne 0 \mbox{ or } b \ne 0 \cr
H^*(\Sigma _{1/2}) \otimes \Lambda ^*H^1(\Sigma _{g-1}) \quad \quad & a=b=0
\end{array}
\right. 
\end{align}
with an appropriate shift in grading on the right-hand side. In both
cases the filtration is equal to the grading modulo a shift. When
$a\ne 0$ or $b\ne 0$ the spectral sequence on $\widehat{CF}$ induced
by $\widehat{CFK}$ collapses at the $E^1$ level and it collapses at
the $E^2$ level when $a=b=0$.
\end{lemma}
Let $Y(g,d)$ the usual groups (see section \ref{notationshort})
defined from $\widehat{HFK} (Y,K,\mathfrak{t}_{a,b})$. In the case
$a=b=0$, the $d_1$ differential from lemma \ref{bmn2} induces a
differential $D$ on $Y(g,d)$, while if $a\neq 0 $ or $b\neq 0$ all
differentials in the spectral sequence starting with $Y(g,d)$ are
trivial since in that case $HFK^\infty\otimes \zee[U,U^{-1}]$ is
isomorphic to $HF^\infty$ of the underlying manifold. To give an explicit
formula for $D$, we decompose $H^* (\Sigma _{1/2})$ as $H^+ (\Sigma
_{1/2}) \oplus H^- (\Sigma _{1/2})$ with
\begin{align} \nonumber
H^+ (\Sigma _{1/2}) = & H^0 (\Sigma _{1/2}) \oplus H^2 (\Sigma _{1/2}) \cr
H^- (\Sigma _{1/2}) = & H^1 (\Sigma _{1/2})
\end{align}
This induces a corresponding decomposition of $Y(g,k)=Y_+(g,k) \oplus
Y_-(g,k) $. One can then check that the differential $D$ becomes
\begin{equation} \label{mnpositive}
D|_{Y_+(g,k)} = 0 \quad \quad D(c\otimes \omega \otimes U^i ) = \omega \otimes U^i + S \wedge \omega \otimes U^{i+1} 
\end{equation}
when $m,n>0$ and it is 
\begin{equation} \label{mnnegative}
D|_{Y_-(g,k)} = 0 \quad \quad D(\omega \otimes U^i ) = c\otimes \omega \otimes U^{i+1} \quad \quad  
D(S\otimes \omega \otimes U^i ) = c \wedge \omega \otimes U^i 
\end{equation}
when $m,n<0$. As in the case of a single twist, there can be no
further differentials in the spectral sequence for $H_*(C\{i<0
\mbox{ and } j\geq k\})$. Assembling all the above we arrive at:
\begin{theorem} \label{main3a}
Let $M(t_\gamma^mt_\delta^n)$ be the mapping torus associated to the
diffeomorphism given by $n$ Dehn twists along $\delta$ followed by $m$
Dehn twists along $\gamma$. We assume that $m\cdot n>0$ and that $\gamma,
\delta \subset \Sigma _g$ are geometrically dual curves supported in a
genus $1$ summand of $\Sigma _g$. Then for $k\ne 0$ the Heegaard
Floer homology of $M(t_\gamma^mt_\delta^n)$ is given by
$$HF^+(M(t_\gamma^mt_\delta^n),\mathfrak{t}_{a,b,k}) = \left\{
\begin{array}{ll} 
 X(g-1,d-1) & a\ne 0 \mbox{ or } b\ne 0  \cr
H_*( Y(g,d), D)  & a=b=0 
\end{array}
\right. 
$$
with $D$ given by \eqref{mnpositive} when $m,n>0$ or
\eqref{mnnegative} when $m,n<0$, and where $d = g-1-|k|$. If in addition
 $g<3k+2$, the above isomorphisms are $\Lambda ^*
 H_1(Y)\otimes _\mathbb{Z} \mathbb{Z}[U]$ module isomorphisms.
\end{theorem}
The homology groups $H_*(Y(g,d),D)$ can be calculated explicitly
using calculations similar to those in section \ref{notationshort}. We
focus on the case $m,n>0$ first. Notice that if we let $c$ denote
a generator of $H^1(\Sigma_{1/2})$ and $S$ a generator of
$H^2(\Sigma_{1/2})$, then one can rewrite the groups $Y_\pm (g,d)$ as
\begin{eqnarray*} Y_+(g,d) &= & X(g-1,d-2) \oplus \left( S \otimes
X(g-1, d) \right) \cr Y_-(g,d)& = & c\otimes X(g-1, d-1 )
\end{eqnarray*}
From this it is easy to see that $Y_-(g,d)$ has no cycles while the
boundaries in $Y_+(g,d)$ are of the form $\omega \otimes U^i +
S\otimes \omega \otimes U^{i+1}$ for $\omega \otimes U^i \in
X(g-1,d-2)$. The homology we seek is then
$$Y_+(g,d) /\{ \omega \otimes U^i + S\otimes \omega \otimes U^{i+1} | \, \omega\otimes U^i  \in X(g-1,d-2)  \}$$ 
The above can easily be seen using arguments as in section
\ref{notationshort} to be isomorphic to $X(g-1,d-1)\oplus \Lambda ^{2g-d-2}
H^1(\Sigma _{g-1})$.

In the case when $m,n<0$, all elements of $Y_-(g,d)$ are boundaries
and don't contribute to the homology. The cycles in $Y_+(g,d)$ are the
elements of the form $\omega \otimes U^i - S \otimes \omega \otimes
U^{i+1}$ (contributing an $X(g-1,d-2)$ summand to the homology) as
well as elements in the \lq\lq bottom-most row\rq\rq in $Y_+(g,d)$,
that is elements of the form $\omega \otimes U^{j-k-1}$ with $\omega
\in \Lambda ^{g+j}H^1(\Sigma _{g-1})$, $j=k,\ldots,g-1$. It is an
explicit check to see that the $X(g-1, d-2)$ summand together with the
elements of the form $\omega \otimes U^{j-k-1}$, $\omega \in \Lambda
^{g+j}H^1(\Sigma _{g-1})$, $j=k+1,\ldots,g-1$ add up to $X(g-1,d-1)$. A
little extra care with gradings leads to
\begin{corollary} With the assumption as in theorem \ref{main3a}, 
the Floer homology of $M(t_\gamma^mt_\delta^n)$ for $m\cdot n >0$ in the
spin$^c$-structures $\mathfrak{t}_{0,0,k}$ (with $k\ne 0$) is
\begin{align} \nonumber
HF^+(M(t_\gamma^mt_\delta^n) , \mathfrak{t}_{0,0,k}) = \left\{
\begin{array}{l}
X(g-1,d-1)\oplus \Lambda ^{2g-d-2}H^1(\Sigma _{g-1})_{( g-1-d)}  \cr
X(g-1,d-1)[-2] \oplus \Lambda ^{2g-d-2}H^1(\Sigma _{g-1})_{( g-1-d)}  
\end{array}
\right.
\end{align}
The first line on the right-hand side above corresponds to the case $m,n>0$ while the second 
line describes the case $m,n<0$. 
\end{corollary}

\subsection{The case  $m \cdot n <0$}

We content ourselves by only
pointing out the major differences that occur here compared with the
previous case.

We may assume that $n>0$. Then the long exact sequence \eqref{LES2} still
applies with the values of the first two terms appropriately adjusted.
An easy case by case analysis (depending on $a$ and $b$) leads to:
\begin{lemma}
Suppose $a$ and $b$ are not both zero. Then the knot Floer homology of
$B(m,n)$ with $m\cdot n<0$ in the spin$^c$-structure
$\mathfrak{t}_{a,b}$ is given by
\[
\widehat{HFK}(L(m,1)\#L(n,1), B(m,n), \T_{a,b}) \cong \zee,
\]
supported in filtration level zero and degree $\ell_{m,a}+\ell_{n,b}$.

If $a = b = 0$, then
\begin{equation} \nonumber
\widehat{HFK}(L(m,1)\# L(n,1), B(m,n), \mathfrak{t}_{0,0},j ) \cong
\left\{
\begin{array}{ll} 
\mathbb{Z}_{( \ell_{n,0} - \ell_{m,0} + 1 )} \quad &
j=1 \cr \mathbb{Z}^3_{( \ell_{n,0} - \ell_{m,0} )} \quad &
j=0 \cr \mathbb{Z}_{( \ell_{n,0} - \ell_{m,0} - 1 )} \quad &
j=-1 \cr 0 & \mbox{otherwise} \end{array} \right.
\end{equation}
where $\ell_{m,a}$ and $\ell_{n,b}$ are defined by \eqref{ella}. The
induced spectral sequence on $\widehat{CF}$ collapses at the $E^1$
stage if either $a\ne 0$ or $b\ne0$ and it collapses at the $E^2$
stage when $a=b=0$.
\end{lemma}
When $m\cdot n<0$, the symmetries of $\widehat{HFK}$ show that the
knot Floer homologies of $B(m,n)$ and $B(-m,-n)$ are equal.

In the following, we will drop the reference to the absolute grading
given above since only the relative grading is relevant in the end
result. However, it will be convenient to suppose that the group
$\widehat{HFK}(B(m,n),j)$ is supported in degree $j$ since with this
convention the gradings on $Y(g,d)$ and $X(g,d)$ are compatible in the
calculations to follow.

In analogy with the previous section, we introduce the genus $3/2$
surface $\Sigma _{3/2}$ whose cohomology groups are $H^k(\Sigma
_{3/2};\mathbb{Z}) = \mathbb{Z}$ for $ k=0,2$ and $H^1(\Sigma _{3/2})
= \mathbb{Z}^3$.
\begin{lemma} \label{bmn3}
For $m\cdot n<0$ the knot Floer homology of $K$ (given by
\eqref{knotdef3}) with respect to the spin$^c$-structure
$\mathfrak{t}_{a,b}$ is
\begin{align} \nonumber
\widehat{HFK}& (Y,K,\mathfrak{t}_{a,b}) = \left\{ 
\begin{array}{ll}
\Lambda ^*H^1(\Sigma _{g-1}) \quad \quad & a\ne 0 \mbox{ or } b \ne 0 \cr
H^*(\Sigma _{3/2}) \otimes \Lambda ^*H^1(\Sigma _{g-1}) \quad \quad & a=b=0
\end{array}
\right. 
\end{align}
with an appropriate shift in grading on the right-hand side. In both
cases the filtration is equal to the grading modulo a shift. When
$a\ne 0$ or $b\ne 0$ the spectral sequence on $\widehat{CF}$ induced
by $\widehat{CFK}$ collapses after the $E^1$ level and it collapses
after the $E^2$ level when $a=b=0$.
\end{lemma}

When $a=b=0$, the $d_1$ differential in $\widehat{HFK}$ can be
described as follows. Since the underlying manifold of $B(m,n)$
has $\widehat{HF}\cong \zee$ and the filtration in $\widehat{HFK}$ is
equal to the grading, the $d_1$ differential must correspond to a
differential on $H^*(\Sigma_{3/2})$ whose homology is $\zee$ supported
in the middle dimension. Let $e$ be a generator for
$H^0(\Sigma_{3/2})$, write $c_1, c_2, c_3$ for the generators of
$H^1(\Sigma _{3/2})$, and write a generator of $H^2(\Sigma _{3/2})$ as
$S$. Then under the $d_1$ differential we can suppose that $S\mapsto c_1$
and that $c_2\mapsto e$ and $c_3\mapsto e$. From this it is easy to
check that for $\omega\in \Lambda^*H^1(\Sigma_{g-1})$, the
differential $D$ on $Y(g,d)$ must be given by:
\begin{align} \label{d1again}
& D(e\otimes\omega \otimes U^i ) = c_1\otimes \omega \otimes U^{i+1} \cr
& D(c_\ell \otimes \omega \otimes U^i ) = \left\{ 
\begin{array}{ll}
0 \quad & \ell =1\cr
e\otimes\omega \otimes U^i - S\otimes \omega \otimes U^{i+1} & \ell =
2,3 \cr
\end{array}
\right. \cr
& D(S\otimes \omega \otimes U^i) = c_1 \otimes \omega \otimes U^i 
\end{align}
There are no further differentials in the spectral sequence for
dimensional reasons, so we get:
\begin{theorem} \label{main3b}
Let $M(t_\gamma^mt_\delta^n)$ be the mapping torus associated to the
diffeomorphism given by $n$ right-handed Dehn twists along $\delta$
followed by $m$ right-handed Dehn twists along $\gamma$. Assume that
$m\cdot n<0$ and that $\gamma, \delta \subset \Sigma _g$ are
geometrically dual curves supported in a genus $1$ summand of $\Sigma
_g$. Then the Heegaard Floer homology of $M(t_\gamma^mt_\delta^n)$ for
$k\ne 0$ is given by
$$HF^+(M(t_\gamma^mt_\delta^n),\mathfrak{t}_{a,b,k}) = \left\{
\begin{array}{ll} 
 X(g-1,d-1) & a\ne 0 \mbox{ or } b\ne 0  \cr
H_*( Y(g,d), D)  & a=b=0 
\end{array}
\right. 
$$
with $D$ given by \eqref{d1again}.   If in addition
 $g<3k+2$, the above isomorphisms are $\Lambda ^*
 H_1(Y)\otimes _\mathbb{Z} \mathbb{Z}[U]$ module isomorphisms.
\end{theorem}
As in the previous section, the group $H_*( Y(g,d), D) $ can be
calculated explicitly. Once again, our first step is to rewrite
$Y(g,d)$ appropriately. We write
\[
H^*(\Sigma_{3/2}) = A_1 \oplus A_0 \oplus A_{-1},
\]
where $A_1 = span\{c_2,c_3\}$, $A_0 = span\{S,e\}$, and $A_{-1}
=span\{c_1\}$. This induces a corresponding splitting of
$\widehat{HFK}(Y,K)$, and thence a decomposition
\[
Y(g,d) = X_1 \oplus X_0 \oplus X_{-1}.
\]
From \eqref{d1again}, we see that $D$ maps $X_i$ into $X_{i-1}$ (and
vanishes on $X_{-1}$). Clearly $D$ maps onto $X_{-1}$. 

The group of cycles in $X_{0}$ is generated by elements of the form
$\xi = e\otimes\omega\otimes U^{i-1} - S\otimes \omega\otimes U^{i}$,
where $\omega\in\Lambda^jH^1(\Sigma_{g-1})$. Such an element is the
boundary of $c_2\otimes\omega\otimes U^{i-1}$ except in the case $i =
0$. This is the situation in which $\xi = S\otimes \omega$ is in
the ``right-hand column'' of $Y(g,d)$; the only such elements that are
cycles are in the bottom right corner. The set of such $\xi$ is
isomorphic to $\Lambda^{2g-2-d}H^1(\Sigma_{g-1})$, and because $S$ is
considered to have degree 1, the contribution to the homology from
$X_0$ is
\[
\Lambda^{2g-2-d}H^1(\Sigma_{g-1})_{(g-1-d)}.
\]
Turning to $X_1$, one checks that the group of cycles is generated by
elements of the form $c_2\otimes\omega\otimes U^i - c_3\otimes
\omega\otimes U^i$. (The only possible case in which more cycles might
appear is in the bottom left corner, but that group contains no
elements of $X_1$.) There are no boundaries, and one can check that
the group spanned by elements of the given form is isomorphic to
$X(g-1,d-1)[-1]$. Together with the previous work, this gives:

\begin{corollary}
The Floer homology of $M(t_\gamma^mt_\delta^n)$ with $m\cdot n <0$ in the
spin$^c$-structures $\mathfrak{t}_{0,0,k}$, with $k\ne 0$ is given as
a group by
\begin{align} \nonumber
HF^+(M(t_\gamma^mt_\delta^n) , \mathfrak{t}_{0,0,k}) = X(g-1, d-1)[-1] \oplus \Lambda
^{2g-d-2} H^1(\Sigma _{g-1})_{( g-1-d)}
\end{align}
where $d = g-1-|k|$.
\end{corollary}

We conclude this section by pointing out that our discussions have focused on the case $k>0$. The fact that the 
results of theorems \ref{main3a} and \ref{main3b} only depend on $|k|$ follows from 
equation \eqref{negativeconversion} together with 
the observation that $\mathfrak{t}_{a,b,-k} = \overline{\mathfrak{t}_{m-a,n-b,k}}$ (see also the comment after
theorem \ref{main2}). 

\section{Dehn twists along a separating curve} \label{four4}
We now turn to the case of a mapping torus $M(t_\sigma^{\pm 1})$ obtained from a
diffeomorphism of $\Sigma_g$ induced by a 
right-handed Dehn twists around a curve $\sigma\subset\Sigma_g$ that
separates $\Sigma_g$ into two components of genera $g_1 = 1$ and $g_2
\geq 1$. For convenience, we focus on the case of a right-handed Dehn twist first; the case 
of a left-handed twist is entirely analogous.
\begin{figure}[htb!] \small
\centering
\includegraphics[width=12cm]{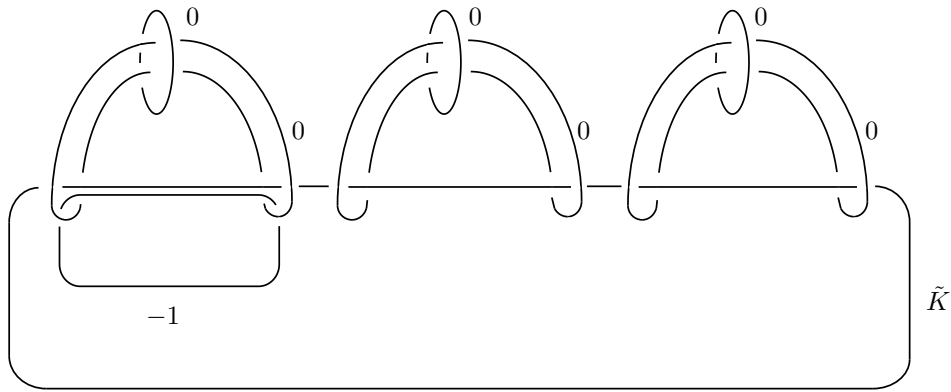}
\put(5,30){$\tilde{K}$}
\put(-60,138){$0$}
\put(-168,138){$0$}
\put(-275,138){$0$}
\put(-18,95){$0$}
\put(-127,95){$0$}
\put(-235,95){$0$}
\put(-290,25){$-1$}
\caption{The knot $\tilde K$ in the case $g=3$.}  \label{sepknot2}
\end{figure}
The surgery picture for $M(t_\sigma)$ is obtained from the usual picture for
$\Sigma_g\times S^1$ by adding the curve
$\sigma$  with framing $-1$. Figure \ref{sepknot2} depicts the case $g
= 3$.

\begin{figure}[htb!] \small
\centering
\includegraphics[width=4cm]{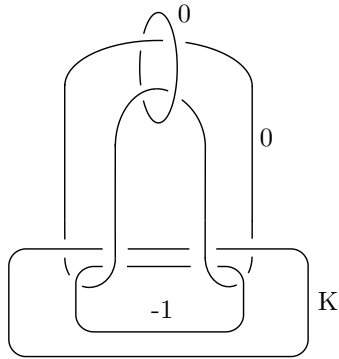} 
\put(-50,127){0}
\put(-19,80){0}
\put(3,18){K}
\put(-60,15){-1}
\caption{The knot $K$ inside  $M = M\{0,0,-1\}$.} \label{sepknot}
\end{figure}

Thus the mapping torus $M(t_\sigma)$ can be seen as the result of 0-framed 
surgery on a knot $\tilde{K}$ (the \lq\lq long\rq\rq circle in figure
\ref{sepknot2}) given by
\[
\tilde{K} = {K} \#^{ g-1} B(0,0)
\]
where $K$ is shown in figure \ref{sepknot}.

Following our usual procedure, we calculate the knot Floer homology of
$\tilde{K}$. To do so, it it follows from 
the connect sum theorem for $\widehat{HFK}$ 
that we need only calculate $\widehat{HFK}(M, K)$, where the 3-manifold 
$M$ is the result of performing $-1$ surgery on one component of the Borromean rings and $0$-surgery on the other
two components.

The homology $\widehat{HF}(M)$ has been calculated in \cite{peter2} (where our $M$ has been   
denoted by $M\{0,0,-1\}$): $\widehat{HF}(M) = \mathbb{Z}^2_{(0)} \oplus \mathbb{Z}^2_{(1)}$. 
We use these groups as in input for calculating the knot Floer homology for $(M, K)$.

\begin{lemma} The knot Floer homology groups for $(M, K)$ are given 
by
\[
\widehat{HFK}(M,K; j) = \left\{\begin{array}{ll} \zee_{(1)} & j = 
1\\ \zee^3_{(0)}\oplus \zee_{(1)} & j = 0 \\ \zee_{(-1)} & j = -1 \\ 
0 & \mbox{otherwise}\end{array}\right.
\]
The spectral sequence that calculates $\widehat{HF}(M)$ from 
$\widehat{HFK}(M,K)$ collapses at the $E_2$ level, and the only 
nontrivial differential is a surjection $d_1: \zee^3_{(0)}\to\zee_{(-1)}$.
\end{lemma}

\begin{proof} We adopt the notation from \cite{ron} where 
$M\{p,q,r\}$ is the manifold gotten by surgery on the Borromean rings
with surgery coefficients $p,q,r$. Notice that $M=M\{0,0,-1\}$.

We write the sequence in knot Floer homology 
arising from the triple 
\[
(M\{0,0,\infty\}, B(0,0))\to (M\{0,0,-1\}, K)\to (M\{0,0,0\}, U)
\]
where $U$ denotes the unknot in $M\{0,0,0\} = T^3$ (since if we replace the $-1$ 
circle in $M$ by a $0$-circle, the knot $K$ slides away from the 
rest of the diagram). This sequence appears as:
\[
\begin{array}{cccc}
\widehat{HFK}(M\{0,0,\infty\}, B(0,0)) & 
\stackrel{F_1}{\to} & 
\widehat{HFK}(M\{0,0,-1\}, K) &
\to\widehat{HFK}(T^3,U)\\
&&&\\
{\begin{array}{ccc} &&\zee_{(1)} \\ & \zee^2_{(0)} & \\ \zee_{(-1)} && \end{array}} & &
{\begin{array}{ccc} &&\zee_{(1)} \\ & B_{(0)} & A_{(1)} \\ \zee_{(-1)} && \end{array} }  &
{\begin{array}{cc} &\\ \quad \zee^3_{(-1/2)} & \zee^3_{(1/2)} \\&\end{array} } 
\end{array}
\]
Here the groups $A$ and $B$ are unknown. The arrays of groups above
are arranged so that the horizontal coordinate corresponds to the
absolute grading of the group and the vertical coordinate is its
filtration level; the maps between arrays preserve the filtration. The
map from the first to the second arrays preserves absolute grading, while
the maps from the second to third and from the third to the first both
decrease absolute grading by 1/2.

Since for any $(Y,K)$ the knot homology $\widehat{HFK}(Y,K)$ is the $E_1$ term in a spectral sequence
associated to a filtration of $\widehat{CF}(Y)$, there is a boundary
map $d_1$ in each of the above arrays mapping one step diagonally down and
to the left. Since $d_1$ is induced from the ordinary boundary in
${CF}$, the homomorphisms in the long exact sequences are
chain maps for $d_1$ (they are maps of spectral sequences).
Furthermore, since when we forget the filtration we have
$\widehat{HFK}(M\{0,0,\infty\}, B(0,0)) = \widehat{HF}(M\{0,0,\infty\})$
and $\widehat{HFK}(T^3, U) = \widehat{HF}(T^3)$, the
$d_1$ boundaries must be trivial in those groups.

The differentials $d_i$, $i\geq 2$ in the middle term must be trivial 
for dimensional reasons. The homology of $d_1$ must
therefore give $\widehat{HF}(M\{0,0,-1\})$, and it follows from the
calculation quoted above that this is $\zee^2_{(1)}\oplus
\zee^2_{(0)}$. In particular, $d_1: B\to \zee$ must be surjective in
order that the homology vanish in degree $-1$; therefore the group $A$ can
have rank 1 or 2 and has rank 1 if and only if $d_1:\zee\to B$ is
trivial. We claim that this boundary map is in fact trivial.

To see this, note that $F_1$ is an isomorphism in filtration level 1.
Furthermore, $d_1 = 0$ in $\widehat{HFK}(M\{0,0,\infty\},B(0,0))$, so
that $F_1\circ d_1 = d_1\circ F_1 = 0$. These two facts combine to 
imply that $d_1:\zee\to B$ must vanish. 
The lemma follows.
\end{proof}
It will be convenient in what follows to write $\widehat{HFK}(M, 
K)$ as
\[
\widehat{HFK}(M,K) \cong \Lambda^*H^1(\Sigma_1)\oplus H^*( S^1).
\]
In this notation, we can express the single nontrivial differential 
in the spectral sequence for $\widehat{HF}(M)$ as the map 
$\Lambda^1H^1(\Sigma_1)\to \Lambda^0H^1(\Sigma_1)$ given by contraction with 
a generator $\gamma$ of $H_1(\Sigma_1)$, which we represent as an embedded 
circle in the torus also denoted $\gamma$.

The connected sum theorem for $\widehat{HFK}$ then gives:

\begin{proposition} The knot Floer homology of $\tilde{K}\subset Y =  M 
\#^{2g-2} (S^1\times S^2)$ is given by
\begin{equation}
\widehat{HFK}(Y,\tilde{K}) = \Lambda^*H^1(\Sigma_{g})\oplus \left[ 
\Lambda^*H^1(\Sigma_{g-1})\otimes H^*( S^1)\right].
\label{septwistHFK}
\end{equation}
\label{septwistHFKprop}
\end{proposition}

Note that in the notation of definition \ref{maxmin}, we have in this
case that $M - m = 1$, where in all previous cases we had $M = m$.

Under the splitting $\Sigma_g = \Sigma_1\# \Sigma_{g-1}$ and
corresponding decomposition
\[
\Lambda^*H^1(\Sigma_g) = 
\Lambda^*_+H^1(\Sigma_g)\oplus\Lambda^*_-H^1(\Sigma_g)
\]
as in section \ref{notationshort}, we have that the $d_1$ differential in 
the spectral sequence for $\widehat{HF}(Y)$ induced by $\tilde{K}$ is 
given by contraction with $\gamma$ on the factor 
$\Lambda^*_-H^1(\Sigma_g)$ in the first summand of 
(\ref{septwistHFK}), and trivial on all other factors. It is not hard 
to see that the homology of $d_1$ is equal to $\widehat{HF}(Y)$, so 
the spectral sequence must collapse at the $E_2$ stage.

We return with the above information to the calculation of 
$HF^+(M(t_\sigma))$. It is a simple matter to check using the results
above that $HF^+_{red}(Y) = HF^+_{red}(Y_p) = 0$ for $p>0$
sufficiently large (recall that $Y = M\# (2g-2)S^1\times S^2$), so the
formalism from section \ref{prelim} applies. As in previous
calculations, we have that $HF^+(M(t_\sigma);k)$ is computed by a
spectral sequence whose $E_2$ term is the homology of $Y(g,d)$ (see
\eqref{three1term}) with respect to a differential constructed from
the differential on $\widehat{HFK}(Y,\tilde{K})$. Here as usual, $d =
g-1-|k|$ and $k\neq 0$. To describe the differential on $Y(g,d)$, note
that \eqref{three1term} together with (\ref{septwistHFK}) give
\[
Y(g,d) = X(g,d) \oplus \left[ X(g-1,d-1)\otimes H^*( 
S^1)\right].
\]
This, together with the results of proposition \ref{septwistHFKprop} 
and the discussion on the differential in 
$Y(g,d)$ from section \ref{notationshort}, easily show that the the 
differential $D$ on $Y(g,d)$ is given by $D_\gamma^+ $ (see \eqref{Dgamma}) on the factor 
$X(g,d)$ and is trivial elsewhere. It follows immediately that 
\begin{equation}
H_*(Y(g,d), D) = H_*(X(g,d), D_\gamma^+ ) \oplus \left[ X(g-1, d-1)\otimes 
H^*( S^1)\right].
\label{septwistanswer}
\end{equation}

A simple argument shows that there can be no higher differentials 
in the spectral sequence calculating $H_*(C\{i<0 \mbox{ and } j\geq 
k\})$, so according to Theorem \ref{morerelevantone}
$HF^+(M(t_\sigma);k)$ is given by \eqref{septwistanswer}.
The homology $H_*(X(g,d), D_\gamma^+)$ was calculated in section
\ref{notationshort}, and this with  similar calculation in the case
$n<0$ leads to:

\begin{theorem} \label{main1}
The Floer homology $HF^+(M(t_\sigma^n) , \mathfrak{s}_k)$ of the mapping
torus $M(t_\sigma^n)$ of a right-handed ($n=1$) or left-handed ($n= - 1$)  Dehn
twist around a genus 1 separating curve $\sigma$  in the surface $\Sigma_{g}$
($g\geq 2$) in the \spinc structure $\mathfrak{s}_k$ whose first Chern
class evaluates to $2k$ on $[\Sigma_g]$ and is trivial on all classes
represented by tori is given by the right hand side of
(\ref{septwistanswer}). More explicitly,
\begin{align} \nonumber
HF^+(M(t_\sigma^n) , & \mathfrak{s} _k) =  \cr 
& (X(g-1,d-1) \otimes H^*(S^1 \sqcup S^1 )){[\varepsilon(n)]} \oplus
\Lambda^{2g-2-d}H^1(\Sigma_{g-1})_{(g-d)}
\end{align}
with $\varepsilon (1) = 0$ and $\varepsilon(-1) =
-1$. Here we assume $k \neq 0$ and set $d = g-1-|k|$.  If in addition
 $g<3k+1$, the above isomorphisms are $\Lambda ^*
 H_1(Y)\otimes _\mathbb{Z} \mathbb{Z}[U]$ module isomorphisms.\qed
\end{theorem}

\Addresses\recd

\end{document}